\documentclass[10pt]{article}
\usepackage[english]{babel}
\usepackage{amsmath}
\usepackage{amssymb}

\newcommand{\Su}{\mathbb{S}^1}
\newcommand{\A}{\mathbb{A}}
\newcommand{\R}{\mathbb{R}}
\newcommand{\N}{\mathbb{N}}
\newcommand{\Z}{\mathbb{Z}}
\newcommand{\Q}{\mathbb{Q}}
\newcommand{\D}{\mathbb{D}}

\newcommand{\C}{{\mathbb C}}
\newcommand{\dimo}{{\bf Proof.}}
\newtheorem{teorema}{Theorem}
\newtheorem{prop}{Proposition}
\newtheorem{cor}[prop]{Corollary}
\newtheorem{lemma}[prop]{Lemma}

\newcommand{\qed}{\hfill {\bf q.e.d.}\medskip}

\begin{document}
\title{Linearization of germs: regular dependence on the multiplier}
\author{Stefano Marmi\thanks{ Scuola Normale Superiore, Piazza dei Cavalieri 7,
56126 Pisa ({\tt s.marmi@sns.it}) } and Carlo Carminati\thanks{
Dipartimento  di Matematica, Universita' di Pisa, Largo Bruno Pontecorvo 5, 56127 Pisa ({\tt carminat@dm.unipi.it})}}
\maketitle
\abstract{ We prove that the linearization of a germ of
holomorphic map of the type
$F_\lambda(z)=\lambda(z+O(z^2))$
 has a ${\cal C}^1$--holomorphic dependence on the multiplier $\lambda$.
${\cal C}^1$--holomorphic functions are
${\cal C}^1$--Whitney smooth functions, defined on compact subsets and which belong to  the kernel of the $\bar{\partial}$
operator.
The linearization is analytic for $|\lambda|\not= 1$ and  the unit
circle~$\Su$ appears as a natural boundary (because of resonances,
{\sl i.e.} roots of unity). However the linearization is still
defined at most points of~${\mathbb S}^1$, namely those points
which lie ``far enough from resonances'', i.e.\ when the
multiplier satisfies a suitable arithmetical condition. We
construct an increasing sequence of compacts which avoid
resonances and prove that the linearization belongs to the
associated spaces of ${\cal C}^1$--holomorphic functions. This is
a special case of Borel's theory of uniform monogenic functions
\cite{Bo}, and the corresponding function space is
arcwise-quasianalytic \cite{MS2}. Among the consequences of these
results, we can prove that the linearization admits an asymptotic
expansion w.r.t. the multiplier at all points of the unit circle
verifying the Brjuno condition: in fact the asymptotic expansion
is  of Gevrey type at diophantine points.} 

\section{Introduction}

A germ of holomorphic diffeomorphism of $(\C ,0)$
\begin{equation}\label{germ}
F_\lambda(z)=\lambda (z + \sum_{k=2}^{+\infty}f_k z^k), \ \ \ \ \ \ (\lambda \in \C^*)
\end{equation}
is {\em linearizable} if there exists a holomorphic germ tangent to the identity
$H_\lambda(z)=z+\sum_{2}^{+\infty} h_k(\lambda) z^k$ which conjugates $F_\lambda$ to the rotation $R_\lambda:z \mapsto \lambda z$ namely
\begin{equation}\label{conjugation}
 F_\lambda\circ H_\lambda=H_\lambda\circ R_\lambda.
\end{equation}
The derivative $\lambda$ of
$F_\lambda$ at the fixed point $z=0$ is called the multiplier of
$F_\lambda$.

If $\lambda$ is not a root of unity there exists a unique formal solution to the conjugacy equation with
coefficients $h_k$, $k\ge 2$,
determined by the recurrence relation
\begin{equation}\label{eq:recurrence}
h_k=\frac{1}{\lambda^{k-1}-1}\sum_{j=2}^{k}f_j
\sum_{\begin{array}{c}\epsilon \in (\Z_+)^j\\ |\epsilon|=k
\end{array}}
h_{\epsilon_1}\cdot ... \cdot h_{\epsilon_j},
\end{equation}
where we follow the usual multi-index notation  $|\epsilon|=\sum_{i=1}^j \epsilon_i$.
Note that $h_k
\in \C(\lambda)[f_2,...,f_k]$.
When $|\lambda|\neq 1$ $F_\lambda$ is always linearizable (by the
classical Koenigs-Poincar\'e  theorem); 
nevertheless the classical estimates
on radius of convergence of $H_\lambda$ deteriorate as $|\lambda|
\to 1$. In the elliptic case, i.e. when $\lambda=e^{2\pi i
\alpha}$ and $\alpha \in \R\setminus\Q \;$, the linearization need
not be convergent due to the contribution from small
denominators in \eqref{eq:recurrence}. After the work of Brjuno
\cite{Br} and Yoccoz \cite{Yo1} we know that all holomorphic germs
with multiplier $\lambda=e^{2\pi i \alpha}$ are analytically
linearizable if and only if $\alpha$ verifies the Brjuno condition
${\cal B}(\alpha)<+\infty$, where   ${\cal B}$ is the Brjuno
function (see the next Section for its
definition and properties).

Let us normalize $F_\lambda$ asking that it is defined and univalent on the unit disk
$\D$. Then one can prove directly, using the classical majorant series' method and Davie's Lemma  
(see \cite{Da} \cite{CM}), that
there are positive constants $b_0,$ $c_0$ (that do not depend on $\alpha$) such that
\begin{equation}\label{crexp}
|h_k| \leq c_0  e^{k ({\cal B}(\alpha)+b_0)}
\end{equation}
where $\lambda = e^{2\pi i \alpha}$, $\alpha \in \R$ and ${\cal
B}$ is the Brjuno function\footnote{It is known that there are
different objects that are called ``Brjuno function'';
nevertheless for this estimate is quite irrelevant which one we
choose, since the difference of two Brjuno functions is bounded by
a universal constant, i.e. independent of $\alpha$ (see Section
\ref{bf}).}.

The same estimate \eqref{crexp} (with larger values $b_0$ and $c_0$)  holds 
uniformly with respect to $\lambda'$ in a cone with vertex in $e^{2\pi
 i \alpha}$.  Therefore for any $\varepsilon
>0$  we will be able to define a closed set $C$ such that there exists $\rho>0$ 
such that
\begin{enumerate}
\item[(i)] ${\rm meas}_2(\C \setminus C)\leq \varepsilon$ and ${\rm meas}_1( C \cap \Su)\geq 2\pi 
-\varepsilon$,
\item[(ii)] for each $\lambda \in C$ the linearization $H_\lambda$ is holomorphic and bounded on $\D_\rho=\{z\in \C \ : \ |z|\leq \rho\}$. 
\end{enumerate}
(Here ${\rm meas}_{d}, \ d\in \{1,2\}$ denotes the $d$-dimensional Lebesgue measure). 

The construction of such a set is performed  removing from $\C$ the union of suitably small connected open
neighbourhoods around the roots of unity\footnote{The property (i) can be realized just asking that 
the ``size'' of the neighbourhood
of each root decays sufficently fast when the order of the root
increases.} and its detailed description can be found in Section \ref{domain}; it will be evident from the construction that 
the radius $\rho$ tends to 0 as $\epsilon$ tends to 0. 

Let us point out that the property (ii) above means that a uniform lower bound on the radius of convergence of $H_\lambda$ holds as $\lambda$ varies in $C$, even near the unit circle.
The set $C$ is sort of a ``bridge'' joining the two connected components of the set of
parameter values  considered in the
Koenigs-Poincar\'e  theorem, crossing the unit circle at some values  $\lambda=e^{2\pi i \alpha}$ with $\alpha$  a Brjuno number.

We address the problem of studying the regularity of this map $\lambda
\mapsto H_\lambda$: we will prove {\it global regularity results}
(see Theorem A below) and {\it local regularity results} (Theorem B).

The global regularity results we prove are inspired by the work of
Borel on uniform monogenic functions \cite{Bo}.
Borel extended the notion of holomorphic
function so as to include functions defined on closed subsets of
$\C$. His  uniform monogenic functions (whose precise definition we
recall and recast in modern terminology in Appendix B) can have, in certain situations, analytic
continuation through what is considered as a natural boundary of
analyticity in Weierstrass' theory. One of Borel's goals was
to determine, with the help of Cauchy's formula, sufficiently
general conditions which would have ensured uniqueness of
the monogenic continuation, {\sl i.e.} a quasianalyticity
property (see \cite{Th}, \cite{Wk} for a modern discussion of ths part).

The importance of Borel's monogenic functions in
parameter-dependent small divisor problems was emphasized by
Kolmogorov \cite{Ko}. Arnol'd discussed in detail this issue in his work~\cite{Ar} on the local
linearization problem of analytic diffeomorphisms of the circle (see \cite{Yo2} for
a very nice introduction and for the most complete results on the subject).
Arnol'd complexified the rotation number but he did not prove that the dependence of the
conjugacy on it is monogenic. This point was dealt with by M.
Herman \cite{He} who also reformulated Borel's ideas using the
modern terminology, Whitney's theory \cite{Wh} on differentiability
of functions on closed sets and the theory of uniform algebras
of (analytic) functions defined on closed sets in the complex
plane. It is Herman's point of view which was developed in
\cite{MS1} and which we will summarize in Appendix B, where we recall the formal definition of ${\cal C}^1 $-holomorphic and ${\cal C}^\infty$-holomorphic functions.
Later Risler \cite{Ris} extended considerably part of Herman's work proving various regularity
results under  less restrictive
arithmetical conditions, namely using the Brjuno conditon as in \cite{Yo2}
instead of a more classical diophantine condition. One should also mention that Whitney
smooth dependence on parameters has been established also in the
more general framework of KAM theory by P\"oschel \cite{Po} who
did not however consider neither complex frequencies nor Brjuno
numbers.

In this paper we will extend the results of Herman and Riesler to
the case of germs of holomorphic diffeomorphisms of $(\C ,0)$. Our
proofs will in fact be more elementary since in this case one can
use a direct approach and the majorant series method applies (see,
e.g. \cite{CM}).

Let us point out that, although the linearization problem makes no
sense for $\lambda=0$ , nevertheless the recurrence
\eqref{eq:recurrence} defines a function $H: \lambda \mapsto H_\lambda$ which turns out 
to be well defined and holomorphic at the
origin: if we denote with ${\cal F}_0= z + \sum_{k=2}^{+\infty}f_k z^k$, 
so that $F_\lambda=\lambda {\cal F}_0$, then $H_0$ turns out to be simply the inverse of ${\cal F}_0$: ${\cal F}_0(H_0(z))=z$.   In fact $H$ can even be extended analitically at infinity
just setting $H_\infty ( z)=z$. Therefore we may consider $H$ as
defined on $C\cup \{\infty\}$ which is  a compact subset of ${\mathbb P}^1\C$:
this has an important consequence since it is proved in \cite{MS2} that the space of ${\cal C}^1$-holomorphic 
functions to which $H$ belongs (see Theorem A below) is  {\em arcwise quasianalytic}\footnote{A function space $X$ is said to be
{\em arcwise quasianalytic} iff the only function that belongs to $X$ and vanishes on an arbitrarily short arc is the null function.}. 

Let us state the main results. In what follows we will assume the germ
$F_\lambda$ to be of the form \eqref{germ}, defined and univalent on
the unit disk $\D$. For any $\rho >0$, $\D_\rho:=\{z\in \C \ : \
|z|\leq \rho\}$ and ${\cal H}^\infty(\D_\rho)$ will denote the complex Banach spaceof functions holomorphic and bounded on $\D_\rho$.

\medskip

\noindent{\bf Theorem A} (Global regularity) {\em $ \ \ $
  For any $\varepsilon>0$ there exist $\rho>0$ and two connected closed sets $C^*$ and $C$  such that
\begin{enumerate}
\item[(a)] $C^*\subset C\subset {\mathbb P}^1\C$, ${\rm meas}_2(\C \setminus C^*)\leq \varepsilon$ and ${\rm meas}_1( C \cap \Su)\geq 2\pi -\varepsilon$
\item[(b)] $H \in C^1_{\rm
hol}(C,{\cal H}^\infty(\D_\rho))$
\item[(b$^*$)] $H \in C^\infty_{\rm
hol}(C^*,{\cal H}^\infty(\D_\rho))$
\end{enumerate}}
\medskip

As will be evident in Section 4.1, for any fixed value $\bar\lambda=
e^{2\pi i \bar\alpha}$ on the unit circle with ${\cal B}(\bar\alpha)<+\infty$
we can manage to build $C^*$ so that $\bar \lambda \in C^*$.
Therefore Theorem A proves that, by a suitable chioice of the set $C^*$, one can extend $H$ and all
its derivatives at any Brjuno point on the circle and 
this leads to the existence of asymptotic
expansions for the linearization $H$ at Brjuno points. 

In fact we can prove that this expansion is
quite 
regular at diophantine points:\medskip


\noindent{\bf Theorem B} (Local regularity) {\em $\ \ $
If $\alpha_0$ is a diophantine point\footnote{Let us recall that an irrational number $\alpha_0$ is diophantine with exponent $\tau_0\geq 2$ if and only if there exists $\gamma>0$ such that for all $p/q\in \Q$ one has 
$|\alpha_0-p/q|\geq \gamma q^{-\tau_0}$.}  with exponent $\tau_0\geq 2$ and $\lambda_0=e^{2\pi i \alpha_0}$, there exists $\rho>0$ such that for any  pair of disks
$\Delta^- \subset \D$ and $\Delta^+ \subset \C\setminus\D$ tangent to $\Su$ in $\lambda_0$  the map
$$\overline{\Delta^+\cup \Delta^-} \ni \lambda \mapsto H_\lambda \in {\cal H}^\infty(\D_\rho)$$
has a Gevrey--$\tau_0$ asymptotic expansion in $\lambda_0$ (we refer the reader to the beginning of Section 5 for its precise definition, see especially \eqref{ggrowth}).
} \medskip

We briefly summarize the content of the paper.  In Section 2 we define
the Brjuno series ${\cal B}$ and we prove several properties of its
sublevel sets.  Since ${\cal B}$ is lower semicontinuous it follows
that the complement of any given  sublevel set $\{x\in \R, \ {\cal B}(x) \leq t\}$
 is the countable disjoint union
of open intervals. These intervals are ``centered'' at those rational values for which the
finite version of the Brjuno series is bounded by the value $t$
defining the sublevel set considered. The discussion of Section 2
prepares the ground for the definition of the domain $C$ where the
conjugation $H:\lambda \mapsto H_\lambda$ is ${\cal C}^1$-holomorphic. The proof of the
$C^1_{\rm hol}$--regularity of $H$ is the main result of Section 3
while in Section 4 we shall restrict the domain of $H$ to a suitably
chosen set $C^*\subset C$ to gain $C^\infty_{\rm hol}$--regularity of
the conjugation.  The proof of Theorem A can be easily obtained
gathering the results of Section 3 (Theorem \ref{c1hol}) and Section
\ref{hr} (Theorem \ref{cihol}). In Section \ref{gevrey} we introduce
 other arithmetically defined real Cantor sets which are useful to
establish the Gevrey regularity of $H$ as claimed in Theorem B.

In the two appendices we recall some elementary properties of the continued fraction expansion of a real 
number (Appendix A) and the definition of ${\cal C}^1 $-holomorphic, ${\cal C}^\infty$-holomorphic and monogenic functions (Appendix B). 

\bigskip

\section{Geometry of the sublevel sets of the Brjuno function}\label{bf}
In the following we shall always denote with $\alpha$ an irrational
number; its continued fraction expansion will be denoted
$\alpha=[a_0,a_1,a_2,...,a_N,...]$ where $a_0\in \Z$ and $a_j\in \Z_+, \ j\geq 1$,
are the partial quotients and $p_j/q_j=[a_0,a_1,a_2,...,a_j]$ is the
$j$-th convergent of $\alpha$.
  If $a_0\in \Z$ and $(a_1, ..., a_N)\in (\Z_+)^N$  we denote $I(a_0,...a_N)$ the set of
real number whose continued fraction begins with the string $(a_0,...,a_N)$. This set is in fact an interval and it is
usually called {\em the cylinder associated to the string of symbols} $(a_0,...,a_N)$
 (for more details and classical
results about continued fraction expansions we refer to appendix A).

To construct the domains on which we will prove the regularity of
the conjugation we shall use the 1-periodic function ${\cal B}$ defined by the following Brjuno series
\begin{equation}\label{b0}
{\cal B}(\alpha):=
\sum_{k=0}^{+\infty}\frac{\log a_{k+1}}{q_k} \ \ \ (\alpha \in \R \setminus \Q), \ \ \ \ \
{\cal B}(r)=+\infty \ \mbox{ if } r\in \Q,
\end{equation}
which is closely related to the classical Brjuno series (see \cite{Br})
$${\cal B}_{cl}(\alpha):=
\sum_{k=0}^{+\infty}\frac{\log q_{k+1}}{q_k}.$$
In fact
it is easily seen that
\begin{equation}\label{bbcl} 0\leq {\cal B}_{cl}(\alpha)-{\cal B}(\alpha)\leq \sum_{k=0}^{+\infty}
\frac{\log(2F_k)}{F_k}<+\infty, \ \ \ \  \mbox{ (where } F_k \mbox{ are the Fibonacci numbers).}
\end{equation}
The above inequalities show that the bound \eqref{crexp} on the growth
of the coefficents of the linearization, which holds for the classical
Brjuno function, must be valid for the Brjuno function ${\cal B}$ as
well (possibly choosing a larger value for the universal constant
$b_0$).  We have chosen ${\cal B}$ instead of ${\cal B}_{cl}$ because
it has various nice properties: its global minimum is 0 and is
attained at the golden mean $ \phi_0:=\frac{\sqrt{5}-1}{2}$.  Moreover
the set $\Phi$ of all local minima of ${\cal B}$ is just the set of
preimages of the golden mean relative to the Gauss map ${\cal
G}(x)=\{1/x\}$ (see also Lemma \ref{lowapp} below):
$$\Phi:=\{\phi \in \R \ : \ {\cal G}^{(n)}(\phi)= \phi_0 \ \mbox{ for some }
n\in \N\}.
$$
Therefore to any element $\phi \in \Phi$ corresponds a continued fraction expansion of the form
$$ \phi= [a_0,...,a_N,1,1,1,...]=[a_0,...,a_N+\phi_0]. $$
In the sequel we will always use only ${\cal B}$.
The following lemmata will be useful to give
a neat description of  sub/super-level sets of the function ${\cal B}$.

\begin{lemma}\label{approssimazione}
If $x\in \R \setminus \Q$ and ${\cal B}(x)<+\infty$ then for all $\varepsilon >0$ exists
$\alpha^{\pm}$ such that
\begin{enumerate}
\item[(i)] $\alpha^-<x<\alpha^+$ and $\alpha^+-\alpha^-<\varepsilon$,
\item[(ii)] $|{\cal B}(\alpha^\pm )-{\cal B}(x)|<\varepsilon$.
\end{enumerate}
\end{lemma}

{\bf Proof:} Let $x:=[a_0,...,a_{N-1}, \ a_N, \ a_{N+1},...]$; fix $N$ odd and
such that
$$ \frac{1}{F_{N}} <\varepsilon, \ \ \ \ \  \ \sum_{k=N}^{+\infty}\frac{\log
a_{k+1}}{q_k} <\varepsilon ,$$ and set
$$\alpha^+:=[a_0,...,a_{N-1}, \ 2 a_N+\phi_0], \ \ \ \
\alpha^-:=[a_0,...,a_{N-1}, \ a_N, \ 2 a_{N+1}+\phi_0].$$
It is clear  that $\alpha^\pm\in I(a_0,...,a_{N-1})$ hence by \eqref{cylinders} 
in Appendix A
$$|\alpha^+-\alpha^-|\leq |\frac{p_{N-1}}{q_{N-1}}-
\frac{p_{N-1}+p_{N-2}}{q_{N-1}+q_{N-2}}|=\frac{1}{q_{N-1}(q_{N-1}+q_{N-2})}\leq
\frac{1}{F_N F_{N-1}}.$$

On the other hand
$$
\varepsilon > \frac{\log 2}{F_n} \geq {\cal B}(\alpha^\pm )-{\cal B}(x) \geq - \sum_{k=N}^{+\infty}\frac{\log a_{k+1}}{q_k}
> -\varepsilon .$$
\qed

\medskip
\begin{lemma}\label{lowapp}
If $x:=[a_0,...,a_{N-1}, \ a_N, \ a_{N+1},...]$ is such that
$a_{2n}>1$ for infinitely many $n\in \N$,
 then for all $\varepsilon > 0$ there exists $\alpha^-$ such that
 \begin{enumerate}
\item[(i)] $x-\varepsilon<\alpha^-<x$,
\item[(ii)] ${\cal B}(\alpha^- )<{\cal B}(x)$.
\end{enumerate}
On the other hand if  for infinitely many $n\in \N$
$a_{2n+1}>1$, then for all $\varepsilon > 0$ there exists $\alpha^+$ such that
 \begin{enumerate}
\item[(i)] $x+\varepsilon>\alpha^+>x$,
\item[(ii)] ${\cal B}(\alpha^+ )<{\cal B}(x)$.
\end{enumerate}
\end{lemma}

{\bf Proof:} It is enough to choose $n$ big enough and such that
$a_{2n}>1$ and set
$$\alpha^-:=[a_0,..., \ a_{2n-1}, \ a_{2n}-1+\phi_0],$$
so that ${\cal B}(x)-{\cal B}(\alpha^-)=-\frac{\log(1-a_{2n}^{-1})}{q_{2n-1}}+
\sum_{k=2n}^{+\infty}\frac{\log a_{k+1}}{q_k} >0$. With a similar trick one can
determine $\alpha^+$.\qed

\medskip

\begin{lemma}\label{Bsci}
The function ${\cal B}$ is lower semicontinuous.
\end{lemma}

{\bf Proof:}
Let $A_t:=\{x\in \R \ :\ {\cal B}(x)>t\}$, $(t\geq 0)$,
  denote the $t$-superlevel  of the function ${\cal B}$. To prove that ${\cal B}$ is semicontinuous
 it is enough to show that $ A_t$ is   open  for all $t\geq 0$.
If ${\cal B}(\bar x)>t$ and $\bar x=[a_0,a_1,a_2,...,a_N,...] \in \R\setminus \Q$ then, for some $N\in \N$,
$\sum_{k=0}^{N}\frac{\log a_{k+1}}{q_k}>t$. Hence ${\cal B}(x) >t$ for all
$x\in I(a_0, ... , a_{N+1})$. A simpler argument settles the case $\bar x$ is rational.\qed

\medskip

It is easy to prove that if $V=]\xi^-,\xi^+[$ is an open interval with
irrational endpoints and $|\xi^+-\xi^-|<1$ then there is a {\em unique} rational point $\bar p/\bar q \in
V$ such that $\bar q <q$ for all $p/q\in V\setminus \{\bar p / \bar q
\}$;
we will call the rational point  $\bar p /\bar q$ the {\em pseudocenter} of
the interval $V$.  Since ${\cal B}$ is lower semicontinuous we know
that the complement of each sublevel is the countable union of
disjoint open intervals with irrational endpoints. Each of these
intervals will be labeled by its pseudocenter. 

In order to characterize the set ${\cal Q}_t$ of pseudocenters of the connected components of the complement of the sublevel $\{ {\cal B}(x)\leq t\}$ we introduce a {\em finite Brjuno function}  (simply
denoted by $B_f$) which is defined on {\em finite continued fractions} by the formula
\footnote{The choice of defining $B_f$ on finite continued fractions instead of $\Q$ avoids the ambiguity which arises from the fact that each rational number admits two different continued fraction expansions, see the Remark below the proof of Lemma 4.} 
$$ B_f([a_0, ..., a_n]):= \sum_{k=0}^{n-1}\frac{\log a_{k+1}}{q_k}.$$

The following lemma gives an accurate description of each of the countable
connected components of the $t$-superlevel sets $A_t$ of the Brjuno function
${\cal B}$ and a precise characterization of the set ${\cal Q}_t$.

\begin{lemma}\label{V}
Let $V=]\xi^-,\xi^+[$ be a connected component of $A_t$ and let
$$\xi^\pm:=[a_0, .. , a_{N-1},\  a_N^\pm, \ a_{N+1}^\pm,...], \ \ \ \ N\geq 1,   \ \ \ a_N^+\neq a_N^- .$$

Then
\begin{enumerate}
\item[(i)]
${\cal B}(\xi^\pm)=t$;
\item[(ii)]
$$\begin{array}{lllll}
a_N^+\geq 2,& a_N^-=a_N^+-1,& a_{N+2k}^+=1 \ \forall k\geq 1,& a_{N+2k+1}^-=1 \ \forall k\geq 0,& \mbox{ if N is even};\\
a_N^-\geq 2,& a_N^+=a_N^--1,& a_{N+2k}^-=1 \ \forall k\geq 1,& a_{N+2k+1}^+=1 \ \forall k\geq 0,& \mbox{ if N is odd}.
\end{array}$$
\item[(iii)] 
The pseudocenter $\bar p/ \bar q$ of $V$ satisfies
\begin{eqnarray*}
  \bar p/\bar q=[a_0,...,a_{N-1},\ a_N^+]=[a_0,...,a_{N-1},\ a_N^-,\ 1],
& B_f([a_0,...,a_{N-1},\ a_N^+])\leq t, & \mbox{ if N is even};\\
  \bar p/\bar q=[a_0,...,a_{N-1},\ a_N^+,\ 1]=[a_0,...,a_{N-1},\ a_N^-],
& B_f([a_0,...,a_{N-1},\ a_N^-])\leq t, & \mbox{ if N is odd}.
\end{eqnarray*}

\item[(iv)] If $p/q\in V$ is a convergent of either $\xi^+$ or $\xi^-$
and $p/q\neq \bar p/\bar q$ then
$$p/q=[a_0,...,a_{N-1},\ a_N^\pm, ...,
a_D^\pm,1],$$ (where $D$ is odd if $p/q$ is a convergent of $\xi^+$
and even if it is convergent of $\xi^-$) and $B_f( [a_0,...,a_{N-1},\
a_N^\pm, ..., a_D^\pm,1])\leq t$ but $B_f( [a_0,...,a_{N-1},\ a_N^\pm,
..., a_D^\pm+1])>t .$
\medskip
\item[(v)] If $p/q \in ]\bar p /\bar q , \xi^+[$ and $p/q$ is {\bf
not} a convergent for $\xi^+$ then exists $p'/q'$ convergent of $\xi^+$
such that $p/q<p'/q'<\xi^+$ and $q'<q$. Moreover, the value of the finite Brjuno
function  exceeds $t$ on both the continued fraction expansions of $p/q$.
A similar statement holds in case $p/q \in ] \xi^-, \bar p /\bar q [$, the only difference
 being that this time $p/q>p'/q'>\xi^-$.
\item[(vi)] If $p/q=[b_0, ..., b_N]$, with $b_N>1$, is a rational number
such that $B_f([b_0, ..., b_N])\leq t$, then $p/q$ is the pseudocenter of
the connected component of $A_t$ which contains it.

\end{enumerate}
\end{lemma}

 $\dimo$
In what follows we will only consider the case $N$ is even since
 the case $N$ odd is symmetric.

{\bf(i)} $\xi^\pm \notin A_t$ implies ${\cal B}(\xi^\pm)\leq t$, on the
other hand it cannot happen that ${\cal B}(\xi^\pm)< t$ because otherwise,
by Lemma \ref{approssimazione}, one could find points in $A_t$ on
which the value of ${\cal B}$ is strictly less than $t$ which is absurd.
\smallskip

{\bf(ii)}
 Let $p_{N-1}/q_{N-1}=[a_0,..., a_{N-1}]$ be the last rational which is a
convergent of both $\xi^\pm$ and stays outside the interval
$]\xi^-,\xi^+[$. Then:
\begin{eqnarray}
a_N^+\geq a_N^-+1 ,\\ \label{a+-}
\frac{p_{N}^-}{q_{N}^-}<\xi^-<\frac{p_{N+1}^-}{q_{N+1}^-}\leq\frac{p_{N}^+}{q_{N}^+}<\xi^+<\frac{p_{N-1}}{q_{N-1}}.
\end{eqnarray}
Setting $\phi:=[a_0,...,a_{N-1}\ , a_N^+-1+\phi_0]$,
 since
$\xi^+>\phi$, and ${\cal B} (\phi)< {\cal B}(\xi^+)\leq t$, we get that
$\xi^-\geq \phi$. Hence $a_N^-\geq a_N^+-1$, in fact by \eqref{a+-}
equality holds.

 If, by contradiction, $a^+_{N+2k}\geq 2$ for some $k
\geq 1$, setting $$\phi:=[a_0,...,a_{N-1}, \ a_N^+, ..., a^+_{N+2k}-1+\phi_0]
$$ we would get $\xi^-<\phi<\xi^+$ while ${\cal B}(\phi)< {\cal B}(\xi^+)\leq t$
which is impossible. An analogous argument shows that $a_{N+2k+1}^-=1$ for all $k\geq 1$.
\smallskip

{\bf(iii)} Let $\bar p/\bar q :=[a_0,...,a_{N-1}, a_N^+]$.  If $p/q \in V$ is a rational number and  $p/q \neq \bar p/\bar q$, then either $p/q\in I(a_0,...,a_{N-1}, a_N^+)$ or $p/q\in I(a_0,...,a_{N-1}, \ a_N^+-1, \ 1)$, and  in both cases
 $q> \bar q$.
\smallskip

{\bf(iv)} If $p/q$ is a convergent of $\xi^+$, and $p/q \neq \bar p/ \bar q$, then $p/q\in ]\bar p/ \bar q, \  \xi^+[$ and
$$p/q=[a_0, ...,a_{N-1}, \ a_N^+, \ ...,a^+_{N+2d-1}, \ 1],   \ \ \ \ (d\geq 1).$$
Moreover, $B_f([a_0, ...,a_{N-1}, \ a_N^+, \ ...,a^+_{N+2d-1}, \ 1])\leq t$. On the other hand, since
$$[a_0, ...,a_{N-1}, \ a_N^+, \ ...,a^+_{N+2d-1}+1+\phi_0]
\in ]\bar p/\bar q, \  p/q[\subset A_t,$$
it follows that
$$B_f([a_0, ..,a_{N-1}, \ a_N^+, ...,a^+_{N+2d-1}+ 1])= {\cal B}( [a_0, ..,a_{N-1}^+, \ a_N^+,...,a^+_{N+2d-1}+1+\phi_0]
)> t.$$
 If $p/q$ is a convergent for $\xi^-$ the argument is symmetric.
\smallskip

{\bf(v)} For $p/q \in ]\bar p /\bar q , \xi^+[$, and $p/q$ {\bf not} a convergent for $\xi^+$, let $[a_0, ...,a_{N-1}, \ a_N^+, \ ...,a^+_{N+2d}]$
be the last convergent smaller than $p/q$, thus $p/q= [a_0, ...,a_{N-1}, \ a_N^+, ...,a^+_{N+2d},  b_1, ...,b_H]$ with $H\geq 1$. We claim that the rational
$p'/q'=[a_0, ..,a_{N-1}, \ a_N^+, ...,a^+_{N+2d}, \ a_{N+2d+1}+ 1]=[a_0, ...,a_{N-1}, \ a_N^+, \ ...,a^+_{N+2d}, \ a_{N+2d+1}, \ 1] $ is the convergent
we are looking for. Indeed, by (i) $p'/q'$ is an even convergent and
 by assumption $p/q<p'/q'$ we deduce that $[b_1, ..., b_H]\geq a_{N+2d+1}+1$, hence
$b_1\geq  a_{N+2d+1}+1$.

Therefore
\begin{eqnarray*}
B_f( [a_0, ..,a_{N-1}, \ a_N^+, ..,a^+_{N+2d}, \ b_1, ..,b_H]) \geq& \!\!\!\!\!\! B_f([a_0, ..,a_{N-1}, \ a_N^+, ...,a^+_{N+2d}, \ a_{N+2d+1}+ 1])>t,
\end{eqnarray*}
the last inequality being a consequence of (iii). If $p/q \in ]\xi^-,\bar p /\bar q[$ the proof can be carried over following the same
argument.

{\bf(vi)} This is a straightforward consequence of the previous statement.
 \qed

\medskip

\noindent {\bf Remark:} ${\cal B}: \R/\Z \to [0,+\infty]$ is
surjective.\\ {\bf Remark:} From now on, if $r\in \Q$, by $B_f(r)$
we will mean the finite Brjuno function
 evaluated on the continued fraction expansion of $r$ which does not end with the figure $1$.
By (vi) of Lemma \ref{approssimazione} the set ${\cal Q}_t$ defined by
\begin{equation}\label{eq:psctrs}
{\cal Q}_t:=\{ r \in \Q \ :\ B_f(r)\leq t\}
\end{equation}
is precisely the set of all pseudocenters of the connected components of $A_t$.

\medskip

By this characterization it is clear that if $t\geq t_0$ then ${\cal
Q}_t \supset {\cal Q}_{t_0}$. It is also interesting and important for the sequel to analyze the
process of disintegration of the connected components of $A_t$; more
precisely let $p/q \in {\cal Q}_{t_0}$ and $]\xi^-_t, \xi^+_t[$ be the
connected component of $A_t$ of pseudocenter $p/q$: the function $t
\mapsto \xi^+_t- \xi^-_t$ is decreasing and has jumps exactly at those
values of $t$ which are image under ${\cal B}$ of local minima in $]\xi^-_{t_0}, \xi^+_{t_0}[$.

\medskip

\begin{prop}\label{dis}
If $t\geq 2t_0$ then all convergents of $\xi^\pm_{t_0}$ belong to ${\cal
Q}_t$.
\end{prop}

$\dimo$ Indeed, if $p/q$ is such a convergent (which is not the pseudocenter of $]\xi^-_{t_0}, \xi^+_{t_0}[$), then
$$ p/q = [a_0,...,a_{N-1}, \ a_N^\pm,..., a^\pm_D,1],$$
$$  \mbox{ where } D\geq 1 \ {\rm is } \ \left\{
\begin{array}{l}
\mbox{ odd for } \xi^+ ,\\
\mbox{ even for } \xi^- .
\end{array}
\right. \ \ \
 \ \ \ $$
Note that, since $ \max (a_N^+, a_N^-) \geq 2$, we have that $$t_0\geq  \max (
{\cal B}(\xi^+), {\cal B}(\xi^-))\geq \frac{\log 2}{q_{N-1}} .$$

So, assuming for the sake of simplicity that $p/q$ is a convergent of $\xi^+$
and dropping the superscript $+$, we can readily check that $B_f([a_0, ..., a_N, ... , a_{D}+1]) \leq t$, which implies that if $p/q\in {\cal Q}_t$ :
\begin{eqnarray*}
B_f([a_0, ... , a_{D}+1])&=& \sum_{k=0}^{D-2} \frac{\log a_{k+1}}{q_k}+\frac{\log (a_D+1)}{q_{D-1}}= \sum_{k=0}^{D-1} \frac{\log a_{k+1}}{q_k}+\frac{\log(1+ a_D^{-1})}{q_{D-1}}\\
&\leq&
t_0 + \frac{\log 2}{q_{N-1}} \leq 2t_0\leq t.
\end{eqnarray*} \qed

Let $M>0$ be fixed. Let $V^*$ be a connected component of $A_{M}$ and let $V$ be a connected component of $A_{3M}$
contained in $V^*$. By the previous remarks we have
\begin{eqnarray*}
V^*= ]\alpha^-, \alpha^+[, & \mbox{ with }& {\cal B}(\alpha^\pm)=M,\\
V= ]\zeta^-, \zeta^+[, & \mbox{ with }& {\cal B}(\zeta^\pm)=3M.
\end{eqnarray*}
We now shall establish a lower bound for  the quantities $|\alpha^+ - \zeta^+|$ and $|\alpha^- - \zeta^-|$.
As usual we will carry over the caculations only for the bound on  $|\alpha^+ - \zeta^+|$, the other case being
analogous. Set $\alpha^+:=[a_0,...,a_{N}, ...], \ \zeta^+:=[c_0,...,c_{N-1}, \ c_N, c_{N+1},\ 1,\ c_{N+3},\ 1,...]$
and let $p/q=[c_0,...,c_{N-1}, \ c_N]$ be the pseudocenter of the interval
$]\zeta^-,\zeta^+[$; we will distinguish the following cases:

\noindent
{\bf Case A: p/q is a convergent of $\alpha^+$ as well}
\begin{quote}
In this case $a_k=c_k$ for all $1\leq k \leq N$. Setting $\phi :=[c_0,..,c_{N-1}, \ c_N, 2c_{N+1},\ 1,\ 1,\ 1,\ 1, ...]$ it is immediate to check
that $p/q<\phi<\zeta^+$ and hence ${\cal B}(\phi)>3M$. Since
$${\cal B}(\phi)= \sum_{k=0}^{N-1}\frac{\log c_{k+1}}{q_k}+\frac{\log 2c_{N+1}}{q_N} \leq {\cal B}(\alpha^+)+\frac{\log 2c_{N+1}}{q_N}
\leq M +\frac{\log 2c_{N+1}}{q_N},$$
it follows that $2c_{N+1}>e^{2Mq_N}$ and so 
\begin{equation}\label{eq:upbnd}
  |\zeta^+-p/q|\leq \frac{1}{c_{N+1}q^2}\leq \frac{2 e^{-2Mq}}{q^2}
\end{equation}
Since  ${\cal B}(\alpha^+)\leq M$, we get $q_N^{-1}\log(a_{N+1})\leq M$ and hence $a_{N+1}\leq e^{Mq_N}$. So 
 \begin{equation}\label{eq:lwbnd}
 |\alpha^+-p/q|\geq \frac{1}{2a_{N+1}q^2}\geq \frac{ e^{-Mq}}{2 q^2}
\end{equation}
Using \eqref{eq:upbnd} and \eqref{eq:lwbnd} we gain
$$  \alpha^+ - \zeta^+\geq  \alpha^+-p/q  - (\zeta^+ -p/q) \geq \frac{e^{-Mq}}{2q^2}(1-4e^{-Mq}),$$
and
\begin{equation}
\alpha^+ -\zeta^+\geq \frac{e^{-Mq}}{4q^2} \ \ \mbox{ holds as soon as } \ \ q\geq \frac{\log 8}{M}
\end{equation}
   \end{quote}

\noindent
{\bf Case B: p/q is not a convergent of $\alpha^+$}
\begin{quote}
If $p/q$ is not a convergent of $\alpha^+$ then  there is some convergent $r/s$  of $\alpha^+$
such that $$p/q<r/s<\alpha^+ \ \ \ \mbox{ and } \ \ \ s<q,$$
hence
$$ \alpha^+-\zeta^+ \geq \alpha^+-r/s\geq \frac{e^{-Ms}}{2s^2} \geq \frac{e^{-Mq}}{2q^2} $$
\end{quote}
The same estimates hold also for $|\alpha^- - \zeta^-|$ so, putting together the cases A and B
we gain the following lemma

\begin{lemma}\label{amz}
There exists a  positive  constant $\nu_0$ such that if $V$ is a
connected component of $A_{3M}$, $p/q$ is the pseudocenter of $V$,
$\zeta \in \partial V$  and ${\cal B}(\alpha)\leq M$ one has
\begin{equation}\label{eq:amz}
|\alpha -\zeta| \geq \nu_0 \frac {e^{-Mq}}{q^2}.
\end{equation}
\end{lemma}

$\dimo$ The proof is straightforward: it is sufficient to note
that the results above imply that \eqref{eq:amz} holds with $\nu_0
= 1/4$ for all but finitely many connected components of $A_t$.
Therefore choosing $\nu_0$ sufficiently small we establish that
\eqref{eq:amz} holds with no exceptions.\qed

\section{${\cal C}^1 $-holomorphic and monogenic regularity of the conjugation}

The main result we shall prove in this section is  that $H\in
{\cal C}^1_{hol}(C_{M};{\cal H}^{\infty}(\D_\rho))$, where $C_{M}$ is a set
obtained removing from the complex plane $ \C$ the union of tiny
neighbourhoods of the roots of unity while $\rho>0$ is suitably
chosen. Let us begin describing the ``domain of regularity'' $C_{M}$.

\subsection{Domain of regularity}\label{domain}
Let $\kappa \in ]0,1[$ be fixed; if $V=]\xi^+, \xi^-[$ is an open
interval in $\R/\Z$ we shall call {\em $\kappa$-diamond\footnote{Or,
simply, {\em diamond}, since we shall not play with different values of $\kappa$.} on V}
the set
$$ \Delta:=\{z\in \C/\Z \ :\ \xi^-<\Re (z) < \xi^+ \ \ \mbox{ and } \ \ |\Im (z)|\leq \kappa \min(\Re z-\xi^-, \xi^+-\Re z)\}.$$
Let $M>0$ be fixed. Let ${\cal Q}_M$ be the set of pseudocenters of connected components of the open superlevel $A_M$ of the Brjuno function  and let
$\Delta(M,r)$ be the diamond on the connected component of $A_M$ containing $r$; it is then easy to check that
$$\Omega_M:=\bigcup_{r\in {\cal Q}_M} \Delta(M,r)$$
is an open neighbourhood of $\Q/\Z$ in $\C/\Z$.
Hence
$K_{M}:=  (\C/\Z) \setminus \Omega_M$ is a closed set which does not contain any rational number.
 Moreover it is straightforward to check that  any of diophantine sets $DC(\gamma, \tau):=\{ \alpha \in \R/\Z \ | \
|\alpha - p/q|\geq \gamma/q^{-\tau} \forall p/q \in \Q\}$ is contained in some $K_M$ for $M$ sufficently large. Since for any fixed $\tau>2$
 ${\rm meas}_1(\R/\Z \setminus DC(\gamma, \tau)) \to 0$ as $\gamma \to 0$ it follows that
${\rm meas}_2 (\C/\Z \setminus K_M) \to 0$ as $M \to +\infty$. It could also be proved (see \cite{MS2}) that ${\rm meas} (K_M) >0$ for all $M>0$ and in fact each point $x\in K_M$ either is isolated (in the exceptional case when $x$ is a local minimunm for the Brjuno function) or $x$ is a point of density for  $K_M$.

\medskip

Let ${\rm exp}^\#:\C/\Z \to \C^*$ be defined as ${\rm exp}^\#(\zeta):= \exp(2\pi i
\zeta)$. For $d>0$ we define $S_d:=\{ \zeta \in \C/\Z \ : \ |\Im (\zeta)|\leq d\}$ the
strip of height $2d$ around the real axis and the annulus
${\rm exp}^\#(S_d):=\A_d$.  We point out that the restriction  ${\rm exp}^\#:S_d \to \A_d$ is a
covering map and locally a biholomorphism, hence by compactness there
exists $\eta>1$ (depending on $d$) such that
\begin{equation}\label{eq:distorsione}
 \eta^{-1}|\zeta-\zeta'|\leq |{\rm exp}^\#(\zeta)-{\rm exp}^\#(\zeta')|
\leq \eta |\zeta-\zeta'| \ \ \ \ \ \ \ \ \forall \zeta,\zeta' \in S_d.
\end{equation}
The set $C_{M}:= {\rm exp}^\#(K_{M})\cup \{0, \infty \}$ is the domain on which the conjugation $H$ will be ``regular''.

\begin{prop}\label{stimahlinf}
There exists a universal constant $b_1$ such that
for all values of the multiplier $\lambda \in C_{M}$ the power series expansion
of the conjugation
\begin{equation}\label{eq:conj}
H_\lambda (z)=\sum_{k=0}^{+\infty} h_k(\lambda)z^k
\end{equation}
has radius of convergence at least $e^{-(M+b_1)}>0$. Moreover $\max_{\lambda \in C_M} |h_k(\lambda)|\leq e^{(M+b_1)k}.$
\end{prop}

{\bf Proof:} We point out that for any fixed $d>0$, if $\lambda \notin
\A_d$ then no small divisor occours in the recurrence
\eqref{eq:recurrence} and the thesis is a straightforward consequence
of the classical majorant series method. Nevertheless the estimates
we get depend on $d$ and deteriorate as $\lambda$ approaches the unit
circle. Therefore we only have to check the statement when $\lambda $
is in some annulus around the unit circle. For this reason in the following we
fix $d>0$ and we consider only those values of the paramenter
$\lambda$ which can be written as $\lambda = e^{2\pi i \xi}$, $\xi \in
K_{M}$, $|\Im \xi|\leq d$. We can associate to any such $\xi$ a point $\xi_0\in
K_{M} \cap \R $ in the following way: $\xi_0:= \Re (\xi) $ if $ \Re
(\xi) \notin A_{M}$ while, if $ \Re (\xi)$ belongs to the connected component  $]\xi^-, \xi^+[$ of
$A_{M}$,
 we shall choose $\xi_0$ to be the nearest point to $\Re(\xi)$ among the two values $\xi^+$
or $\xi^-$; we define also $\lambda_0 = e^{2\pi i \xi_0}$.

In this way  we can easily check that
$$|\lambda^k-1|\geq \delta |\lambda_0^k-1|, \ \ \ \ \forall \lambda \in C_{M}, $$
where $\delta:= \min(e^{-d}\eta^{-2}(1+\kappa^{-2})^{-\frac{1}{2}}, \frac{1-e^{-d}}{2})$.
By the recurrence relation \eqref{eq:recurrence} we get that
\begin{equation}\label{confronto}
|h_k(\lambda)|\leq\delta^{-k} |h_k(\lambda_0)|.
\end{equation}

By the Brjuno estimate \eqref{crexp} (see also \cite{CM} for its
 proof) we deduce that for all $\lambda_0
 \in C_M\cap{\mathbb S}^1$ the radius of convergence $\rho(\lambda_0)$ of
 the series $H_{\lambda_0}(z)=\sum_{k=0}^{+\infty} h_k(\lambda_0)z^k$
 satisfies $\rho(\lambda_0)\geq e^{-(M+b_0)} >0$.  So, by \eqref{confronto}, we get that for
 all $\lambda\in C_{M}\cap\A_d$ the radius of convergence of series
 $H_\lambda (z)=\sum_{k=0}^{+\infty} h_k(\lambda)z^k$ is greater than
 $\delta^{-1} e^{-(M+b_0)}$. \qed

From now on we set $d:=\kappa/2$, so that $\Delta(3M, p/q)\subset S_d$
for all $p/q \in \Q$ and $\eta$ will be the constant appearing in
\eqref{eq:distorsione} relative to $S_d$.

\begin{lemma}\label{gap}
Let $p/q \in \Q$, $\Delta(3M, p/q)$ be the diamond on the connected component of the superlevel $A_{3M}$ containing
$p/q$ and let $D(3M, p/q):={\rm exp}^\#(\Delta(3M, p/q))$. Then, if $\lambda \in C_{M}$
$$ d(\lambda, \partial D(3M, p/q)) \geq \nu_1\frac{e^{-Mq}}{q^2},$$
 where $ \nu_1= \frac{\nu_0}{\eta} \frac{\kappa}{\sqrt{1+\kappa^2}}$
\end{lemma}
{\bf Proof:}Let $\lambda = e^{2\pi i \alpha}, \ \alpha \in K_{M}$,  we immediately get
$$ d(\lambda, \partial D(3M, p/q)) \geq \eta^{-1} d(\alpha, \partial \Delta(3M, p/q)) 
\geq \nu_1\frac{e^{-Mq}}{q^2},$$
where the last inequality follows from Lemma \ref{amz} together with an elementary geometrical argument. \qed

\subsection{${\cal C}^1_{hol}$--regularity}
With a slight abuse of notation let us set $H_k:= \max_{\lambda \in C_{3M}}
|h_k(\lambda)|$; we know from Proposition \ref{stimahlinf} that the series
$\sum H_k z^k$ has a positive radius of convergence bounded from below by
$\rho_0:=e^{-(3M+b_1)}$.
The main result in this section is the following:
\begin{teorema}\label{c1hol}
Let $\rho\in ]0, e^{-2M}\rho_0[$. Then the map $h: \lambda \mapsto H_\lambda$
belongs to the space of functions ${\cal C}^1_{hol}(C_{M};{\cal H}^{\infty}(\D_\rho)).$
\end{teorema}
 We already know that, by  virtue of Proposition \ref{stimahlinf},
\begin{equation}\label{eq:normconv}
H_\lambda (z)=\sum_{k=1}^{+\infty} h_k(\lambda) z^k
\end{equation}
has positive radius of convergence for $\lambda \in C_M$, moreover
each of the coefficients $h_k$, defined by the recurrence relation
\eqref{eq:recurrence}, is a rational function in the variable
$\lambda$ and is holomorphic away from the roots of unity of order
strictly less than $k$.

In order to prove the theorem we shall show that the series
\eqref{eq:normconv} is normally convergent in
${\cal C}^1_{hol}(C_{M};{\cal H}^{\infty}(\D_\rho))$; and this will be a
straightforward consequence of point (iii) of the next lemma:

\begin{lemma}
There exists a positive constant $ L >1$ such that
\begin{enumerate}
\item[(i)]
$|h_k'(\lambda)|\leq L (1+k^4) H_k e^{2M k} \ \ \ \forall \lambda \in C_{M}$
\item[(ii)]
$ \displaystyle \left|\frac{h_k(\lambda_1)-h_k(\lambda_0)}{\lambda_1-\lambda_0} -h'(\lambda_0)\right|\leq 2L (1+k^4) H_k e^{2M k} \ \ \ \forall \lambda_1,\ \lambda_0 \in C_{M}$
\item[(iii)]
$\|h_k(\lambda) z^k\|_{{\cal C}^1_{hol}(C_{M};{\cal H}^\infty(\D_\rho))} \leq 4 L (1+ k^4)
  H_k (e^{2M}\rho)^k$.
\end{enumerate}
\end{lemma}

{\bf Proof:}\\
From now on $r$ will always denote a rational number and ${\rm ord}(r):=\min\{n\in \N^* \ : nr\in \Z\}$; moreover we use the following notation:

\begin{equation}\label{eq:qk}
\begin{array}{ll}
\Q_k=\{r\in \Q: \ {\rm ord}(r)< k\} ; &    D_r:= {\rm exp}^\#(\Delta (3M,r)); \\
{\cal Q}_{t,k}:=\{ r \in \Q \ :\ B_f(r)<t \} \cap \Q_k; &
\end{array}
\end{equation}

Since $h_k$ is a rational function with poles located on the roots of unity of order less than $k$
$${\cal R}_k:={\rm exp}^\#(\Q_k) \subset \bigcup_{r\in {\cal Q}_{3M,k}}D_r$$
For any $R>|\lambda|$ we get, by Cauchy formula,
\begin{equation}\label{eq:cauchy}
h_k(\lambda)=\sum_{r\in{\cal Q}_{3M,k}} \frac{1}{2\pi i}\int_{\partial D_r} \frac{h_k(\zeta)}{\zeta-\lambda} d\zeta + \frac{1}{2\pi i}\int_{\partial \D_R}\frac{h_k(\zeta)}{\zeta-\lambda} d\zeta.
\end{equation}
In fact, letting $R \to +\infty$, we realize that the term
$\int_{\partial \D_R} \frac{h_k(\zeta)}{\zeta-\lambda} d\zeta$ must
vanish; hence in the following we will always neglect this term.  We
can wrtite the integral representation both for the derivative
of $h_k$ and for the Taylor remainder ${R}_2(h_k,\lambda_0, \lambda_1):=
h_k(\lambda_1)-h_k(\lambda_0)-h'_k(\lambda_0)(\lambda_1-\lambda_0)$:
\begin{equation}\label{eq:cauchy1}
h'_k(\lambda)=\sum_{s\in{\cal Q}_{3M,k}} \frac{1}{2\pi i}\int_{\partial D_s} \frac{h_k(\zeta)}{(\zeta-\lambda)^2} d\zeta,
\end{equation}
\begin{equation}\label{eq:taylor1}
\left|(\lambda_1-\lambda_0)^{-1}{R}_2(h_k,\lambda_0, \lambda_1)\right|=\sum_{s\in{\cal Q}_{3M,k} } \frac{1}{2\pi i}\int_{\partial D_s}
\left[
\frac{1}{(\zeta-\lambda_1)(\zeta-\lambda_0)} -\frac{1}{(\zeta-\lambda_0)^2}
\right]
h_k(\zeta) d\zeta .
\end{equation}
So we get the estimates

\begin{equation}\label{eq:scauchy1}
|h'_k(\lambda)|\leq\frac{H_k}{2\pi}\sum_{s\in{\cal Q}_{3M,k}}d(\lambda, \partial D_s)^{-2}\int_{D_s}|d\zeta|,
\end{equation}
\begin{equation}\label{eq:staylor1}
\left|(\lambda_1-\lambda_0)^{-1}{R}_2(h_k,\lambda_0, \lambda_1)\right|
\leq \frac{H_k}{2\pi}\sum_{s\in{\cal Q}_{3M,k}}2d(\lambda, \partial D_s)^{-2}\int_{D_s}|d\zeta|.
\end{equation}
Lemma \ref{amz} gives an upper bound on the term $  d(\lambda, \partial D_s)^{-2}$
$$\int_{D_s}|d\zeta|\leq \eta
\int_{\Delta_s}|dz|\leq 2\eta (1+\kappa^{-2})^{1/2}|V_s| \ \ \ \ {\rm for } \  \ \ s\in {\cal Q}_{3M,k}.$$
Since $\sum_{s\in {\cal Q}_{3M,k}} |V_s|<1$ we finally get
\begin{eqnarray*}
|h'_k(\lambda)|&\leq& \frac{H_k}{2\pi}[2\eta(1+\kappa^{-2})^{1/2}\nu_1^{-2}k^4e^{Mk}],\\
\left|(\lambda_1-\lambda_0)^{-1}{R}_2(h_k,\lambda_0, \lambda_1)\right| &\leq& \frac{H_k}{\pi}[2\eta(1+\kappa^{-2})^{1/2}\nu_1^{-2}k^4e^{Mk}].
\end{eqnarray*}
The thesis follows choosing $L:=\max\{1,
 \frac{\eta}{\pi}(1+\kappa^{-2})^{1/2} \nu_1^{-2} \}$. \qed

\subsection{Monogenic regularity.}
We refer the reader to Appendix B for a more detailed treatment of
monogenic functions. 
Let us choose an increasing sequence of positive
values $M_j$ such that $\lim M_j = +\infty$ and set $C_j:= C_{M_j}$.

Consider the Banach space $B_\ell=\cap_{j=0}^\ell {\cal C}^1_{hol} (C_j, {\cal
H}^\infty (\D_{r_j}))$, where $r_j=e^{-(3M_j+b_1}$ with the norm $\Vert f\Vert_{B_\ell} =
\max_{0\le j\le l} \Vert f\Vert_{{\cal C}^1_{hol} (C_j, {\cal
H}^\infty (\D_{r_j}))}$. Clearly the injections $i_\ell\, :\,
B_\ell\hookrightarrow B_{\ell-1}$ are bounded linear operators between
Banach spaces with norms $\Vert i_\ell\Vert\le 1$. The projective
limit of the system of Banach spaces
$$ {\cal M}( (C_\ell),
\C\{z\}) = \varprojlim B_\ell
$$
is a space of monogenic functions with values in the holomorphic
germs $\C\{z\}$. This is a Fr\'echet space with the family of seminorms
$(\Vert \cdot\Vert_{B_\ell})_{\ell\in \N}$. 

Thus Theorem \ref{c1hol} has the following corollary
\begin{cor}
Let $(C_\ell)$ be as above. Then the linearization H belongs to the space
${\cal M}( (C_\ell),
\C\{z\})$
 of $\C\{z\}$-valued monogenic functions.
\end{cor}

\section{Higher regularity}\label{hr}
Using the Cauchy formula \eqref{eq:cauchy} we get an integral representation for the $m$-th derivative of $h_k$ as well:
\begin{equation}\label{eq:cauchym}
h^{(m)}_k(\lambda)=\sum_{s\in{\cal Q}_{3M,k} } \frac{m!}{2\pi i}\int_{\partial D_s} \frac{h_k(\zeta)}{(\zeta-\lambda)^{m+1}} d\zeta .
\end{equation}
It is quite easy to see that, if we just  followed the same lines of the previous section, in order to gain ${\cal C}^m_{hol}$ regularity we would have to shrink
the radius $\rho$ of the disk $\D_\rho$ and trying to prove that $h$ is
$C^\infty$-holomorphic  would lead to a disk of convergence of radius zero.

To avoid this problem we use an idea of \cite{Ris}: we will prove that
$h \in C^\infty_{ hol} (C^*,{\cal H}^{\infty}(\D_\rho))$ for some $\rho>0$, where
this time $C^*$ will be  somewhat smaller than the set $C_{M}$ considered in the
previous section.

\medskip
\subsection{Construction of the domain $C^*$}
Let $(M_n)_{n\in \N}$ be a decreasing sequence such that $M_n \to 0$ as
$n \to \infty$ and define the set
$$K^*_{(M_n)}:=\{ x\in \R/\Z \ :\ \sum_{k=n}^{+\infty} \frac{\log a_{k+1}}{q_k}\leq M_n\ \forall n \in \N \}.$$
Of course $K^*_{(M_n)}\subset K_M$ for $M=M_0$ and it is not difficult to check that $K_{(M_n)}$ is compact in $\R/\Z$.
It is clear that $K_{(M_n)} \neq \emptyset$ and it is possible to choose $M_n\to 0$ such that meas$(K^*_{(M_n)})>0$: in fact for all $\tau, \gamma$ there exists $M_n\to 0$ such that $K^*_{(M_n)}\supset DC(\gamma, \tau)$ (see also Section \ref{gevrey}).
As before, we define $A^*:=(\R/\Z)\setminus K^*_{(M_n)}$ and we denote with ${\cal Q}^*$ the set of all pseudocenters of the connected components of $A^*$. We now have a list of technical lemmata that
will be useful later on.

\begin{lemma}\label{priori}
 There exists a function $ Q :\N \to \R$ such that if
$\alpha :=[a_0, ..., a_n, ...] \in K^*_{(M_n)}, $ and $p_n/q_n= [a_0, ..., a_n]$ is the $n$-th convergent
then $F_n \leq q_n \leq Q_n$ (where $F_n$ are, as usual, the Fibonacci numbers).
\end{lemma}
{\bf Proof} The sequence $Q_n$ defined by the recurrence
$$\left\{
\begin{array}{l}
Q_0=1,\\
Q_{n+1}=e^{M_n Q_n} Q_n + Q_{n-1},
\end{array}
\right.
$$
does the job.\qed
\medskip

Next lemma is almost a clone of Lemma \ref{V}.

\begin{lemma}\label{V*}
Let $V^*=]\alpha^-,\alpha^+[$ be a connected component of $A^*$ and let
$$\alpha^\pm:=[a_0, .. , a_{N-1},\  a_N^\pm, \ a_{N+1}^\pm,...], \ \ \ \ N\geq 1,   \ \ \ a_N^+\neq a_N^-.$$
Then
\begin{enumerate}
\item[(i)]
$$\begin{array}{lllll}
a_N^+\geq 2,& a_N^-=a_N^+-1,& a_{N+2k}^+=1 \ \forall k\geq 1,& a_{N+2k+1}^-=1 \ \forall k\geq 0,& \mbox{ if N is even};\\
a_N^-\geq 2,& a_N^+=a_N^--1,& a_{N+2k}^-=1 \ \forall k\geq 1,& a_{N+2k+1}^+=1 \ \forall k\geq 0,& \mbox{ if N is odd}.
\end{array}$$
\item[(ii)]
There is a unique rational number $\bar p/\bar q$ which is a convergent
  of both $\alpha^\pm$;
 $\bar p/\bar q$ is the rational number with lowest denominator in $V^*$ (it is then the pseudocenter of $V^*$) and
\begin{eqnarray*}
\bar p/\bar q=[a_0,...,a_{N-1},\ a_N^+]=[a_0,...,a_{N-1},\ a_N^-,\ 1], & B_f([a_0,...,a_{N-1},\ a_N^+])\leq M_0, & \mbox{ if N is even};\\
\bar p/\bar q=[a_0,...,a_{N-1},\ a_N^+,\ 1]=[a_0,...,a_{N-1},\ a_N^-], &B_f([a_0,...,a_{N-1},\ a_N^-])\leq M_0, & \mbox{ if N is odd}.
\end{eqnarray*}
\item[(iii)] If $p/q\in V^*$ is a convergent of either $\alpha^+$ or $\alpha^-$
and $p/q\neq \bar p/\bar q$ then
$$p/q=[a_0,...,a_{N-1},\ a_N^\pm, ...,
a_D^\pm,1]$$ (where $D$ is odd if $p/q$ is a convergent of $\alpha^+$
and even if it is convergent of $\alpha^-$).
\medskip
\item[(iv)] If $p/q \in ]\bar p /\bar q , \alpha^+[$ and $p/q$ is {\bf
not} a convergent for $\alpha^+$ then exists $p'/q'$ convergent of $\alpha^+$
such that $p/q<p'/q'<\alpha^+$ and $q'<q$.
A similar statement holds in case $p/q \in ] \alpha^-, \bar p /\bar q [$, the only difference
 being that this time $p/q>p'/q'>\alpha^-$.
\medskip
\item[(v)] Any convergent of $\alpha^\pm$ is a pseudocenter of some connected component of $A_{3M_0}$
(hence belongs to ${\cal Q}_{3M_0}$).
\end{enumerate}
\end{lemma}
{\bf Proof:} We just sketch some details, since the whole proof of
(i)-(iv) is just a repetition the arguments of Lemma \ref{V} while (v)
follows from the same argument as in Proposition \ref{dis}.  Let
$V^*=(\alpha^-,\alpha^+)$ be a connected component of $A^*$ and let
$\bar r$ be the pseudocenter of $V^*$. Write $\bar r=[a_0, ...,a_N]$
with $a_N\geq 2$ and assume, just to fix ideas, that $N$ is
even. Letting $r^-:=[a_0,..., a_{N-1}]$ and $r^+:=[a_0, ... , a_N -1]$,
it is readily checked that ${\rm ord}(r^\pm)< {\rm ord}(\bar r)$ and
$r^-<\bar r<r^+$, hence, by the minimality of the order of $\bar r$,
$V^*\subset (r^-, r^+)$. On the other hand, for any fixed $n\in \N$,
the expression $\sum_{k=n}^{+\infty} \frac{\log a_{k+1}}{q_k}$ attains
its minimum value on the interval $(\bar r, r^+)$ at the point
$\phi^+:=[a_0, ..., a_N+\phi_0]$ while the minimum value on $(
r^-,\bar r)$ is attained at $\phi^-:=[a_0, ..., a_N-1+\phi_0]$.
This implies that $(\alpha^-, \alpha^+)\subset (\phi^-,\phi^+)$ and $\bar r=[a_0,...,a_N]$ 
is a common convergent of both $\alpha^\pm$.

\qed
\medskip

Let $k\in \N$ be fixed, $r \in {\cal Q}^* \cap \Q_k$ (recall the definitions \eqref{eq:qk}) and let $V^*_r=]\alpha^-, \alpha^+[$  be the connected component of $A^*$ with pseudocenter $r$.
Let us define
$$
s^+(k,r):=\max\{s\in \Q_k\cap]\alpha^-,\alpha^+[\}, \ \ \ \ \ s^-(k,r):=\min\{s\in \Q_k\cap]\alpha^-,\alpha^+[\}.
$$
By virtue of (iv) of the previous lemma $s^+$ is a convergent of $\alpha^+$, while $s^-$ is a convergent of $\alpha^-$ therefore, by (v) of the previous lemma,  both $s^\pm$ belong to ${\cal Q}_{3M_0,k}$.
Let
$$
s(k,r):=\left\{
\begin{array}{ll}
s^+(k,r) & \ {\rm if}\ d(\alpha^+, \partial V_{s^+(k,r)})<d(\alpha^-, \partial V_{s^-(k,r)})  \\
s^-(k,r) & \ {\rm if}\ d(\alpha^+, \partial V_{s^+(k,r)})\geq d(\alpha^-, \partial V_{s^-(k,r)})
\end{array}
\right.
\ \ \ \ {\rm and} \ \ q(k,r):= {\rm ord}(s(k,r)) < k.
$$
Note that if $r \in {\cal Q}^* \cap \Q_k$ then ord$(r)\leq q(k,r)<k$

Thus, if $s\in {\cal
Q}_{3M_0,k}\cap V^*_r$ and $V_s$ is the connected component of $A_{3M_0}$ of pseudocenter $s=p/q$, we have that
$$d(\alpha, \partial V_s)\geq d(\alpha, \partial V_{s(k,r)}).$$

\medskip

The first part of next lemma is just  Lemma \ref{amz}, the second contains the extra information
we shall need to prove that $h$ is ${\cal C}^\infty_{hol}$.
\medskip

\begin{lemma}\label{amz*}
\smallskip

\begin{enumerate}
\item
Let $s \in {\cal Q}_{3M_0}$ and $V_s$ be the connected component of pseudocenter $s=p/q$, then
$$ d(\alpha, \partial V_s)\geq \nu_0 \frac{e^{-M_0 q}}{q^2}, \ \ \ \ \forall \alpha \in K^*_{(M_n)} $$
(where $\nu_0$ is a constant independent of $s$).
 \item
Let $n, k\in \N$,  $V^*_r=]\alpha^-,\alpha^+[$ be the connected component of $A^*$ with pseudocenter $r$
and assume that $q(k,r)\geq \max\{Q_n, \frac{\log 8}{M_0} \}$.
If $s\in {\cal
Q}_{3M_0,k}\cap V^*_r$ then
\begin{equation}\label{eq:amz*}
d(\alpha, \partial V_s)\geq \frac{e^{-M_n k}}{4k^2} \ \ \ \  \forall \alpha \in K^*(M_n).
\end{equation}
\end{enumerate}
\end{lemma}

{\bf Proof:}
We just have to prove the second statement; for sake of simplicity let us assume that $s(k,r)=s^+(k,r)$, the other case being
analogous. If $s\in {\cal
Q}_{3M_0,k}\cap V^*_r$  then
$$d(\alpha, \partial V_s)\geq d(\alpha, \partial V_{s(k,r)})\geq d(\alpha^+, \partial V_{s^+(k,r)}) .$$
Setting $V_{s^+(k,r)}:=]\zeta^-, \zeta^+[$, $s^+(k,r):=p/q$ (so that $q=g(k,r)\geq Q_n$)  we can repeat the argument of the end of subsection \ref{domain}:
$$|\alpha^+-\zeta^+|\geq|\alpha^+ -\frac{p}{q}|-|\frac{p}{q}-\zeta^+|\geq \frac{e^{-M_nq}}{2q^2}-2\frac{e^{-2M_0q}}{q^2}\geq
 \frac{e^{-M_nq}}{2q^2}[1-4e^{-M_0q}]\geq \frac{e^{-M_nq}}{4q^2} \geq \frac{e^{-M_n k}}{4k^2}.$$
\qed

\bigskip

If $n \in N$ is fixed and $k \geq \frac{\log 8}{M_0}$, the following decomposition shall be useful
${\cal Q}^*={\cal Q}^*_1(k,n)\cup {\cal Q}^*_0(k,n)$ where
\begin{equation}\label{qs10}
{\cal Q}^*_1(k,n):=\{r \in {\cal Q}^* \ : q(k,r) \geq Q_n \}, \ \ \ \ \
{\cal Q}^*_0(k,n):=\{r \in {\cal Q}^* \ : q(k,r) < Q_n \}
\end{equation}

Let $\kappa \in ]0,1[$ be fixed and let us carry over the construction of
subsection \ref{domain}: if $V^*_r$ is a connected component $A^*$ and $\Delta^*_r$ is the $\kappa$-diamond over $V^*_r$  we  call
$$\Omega^*:=\bigcup_{r\in {\cal Q}^*} \Delta^*_r$$
it is then easy to check that
$\Omega^*$ is an open neighbourhood of $\Q/\Z$.

 The closed set $K^*:=  (\C/\Z)\setminus \Omega^*$ is connected and does
not contain any rational number.
\medskip

The set $C^*:= {\rm exp}^\#(K^*)\cup \{\infty\}$ will be the domain on which the conjugation will be Whitney smooth.

\subsection{Proof of the regularity}

As before, let $H_k:= \max_{\lambda \in C_{3M_0}} |h_k(\lambda)|$, and let
$\rho_0=e^{-(3M_0+b_1)}$ be the radius of convergence of $\sum H_k z^k$.
We are now able to show that $h\in {\cal C}^\infty_{hol}(C^*;{\cal H}^\infty(\D_\rho))$ for every
 $\rho\in ]0,e^{-1}\rho_0[$.

\begin{teorema}\label{cihol}
For all  $n, m \in \N$
there exist constants $L_m, \ \Lambda_{m,n}$ satisfying
\begin{enumerate}
\item $\displaystyle \max_{\lambda \in C^*} |h_k^{(m)}(\lambda)| \leq L_mk^{2m+2} H_k e^k, \ \ \ \ \forall \lambda \in C^*;$
\item $\displaystyle \sup_{\lambda_0, \lambda_1 \in C^*} \frac{|{R}_{n+1} (h_k^{(m)}, \lambda_1,\lambda_0)|}{|\lambda_1-\lambda_0|^n}| \leq 
\Lambda_{m,n}k^{2m+2n+2} H_k e^k,   \ \ \ \ \forall \lambda \in C^*,$
\end{enumerate}
where $${R}_{n+1} (f, \lambda_1, \lambda_0):= f(\lambda_1)- \sum_{j=0}^{n}
 \frac{f^{(j)}}{j!} (\lambda_0)(\lambda_1-\lambda_0)^j,$$
 is the Taylor remainder of order $n+1$.
\end{teorema}

Before plunging into the proof let us remark that the theorem implies 
$h\in{\cal C}^\infty_{hol}(C^*;{\cal H}^\infty(\D_\rho))$ as soon as $\rho\in ]0,e^{-1}\rho_0[$; indeed, (1) and (2)
imply that the series$\sum_{k=0}^{+\infty}h_k(\lambda) z^k$ is normally convergent in ${\cal C}^m_{hol}(C^*;{\cal H}^\infty(\D_\rho))$ for all $m \in \N$.

{\bf Proof} We have already seen that
$$
h^{(m)}_k(\lambda)=\sum_{s\in{\cal Q}_{3M_0,k} }
\frac{m!}{2\pi i}\int_{\partial D_s}
\frac{h_k(\zeta)}{(\zeta-\lambda)^{m+1}} d\zeta,
$$
hence
\begin{equation}\label{eq:hm}
  |h^{(m)}_k(\lambda)|\leq H_k\frac{m!}{2\pi }\sum_{s\in{\cal Q}_{3M_0,k} } \int_{\partial D_s} d(\zeta,C^*)^{-m-1} |d\zeta|.
\end{equation}
Moreover we can write down an explicit expression for the Taylor
remainder of $\phi_\zeta(\lambda):=\frac{1}{(\zeta -\lambda)^{m+1}}$:
$${R}_{n}(\phi_\zeta,\lambda_1, \lambda_0)=
 \sum_{k=0}^m \binom{k+n-1}{k}(\zeta -\lambda_1)^{k-m-1}(\zeta-\lambda_0)^{-n-k}(\lambda_1-\lambda_0)^n,$$
and we shall use the form
\begin{eqnarray*}
  {R}_{n+1}(\phi_\zeta,\lambda_1, \lambda_0)&=&  {R}_{n}(\phi_\zeta,\lambda_1, \lambda_0)- \frac{\phi_\zeta^{(n)}(\lambda_0)}{n!}(\lambda_1-\lambda_0)^n\\
  &=&
  (\lambda_1-\lambda_0)^n \left[ \sum_{k=0}^m \binom{k+n-1}{k}(\zeta -\lambda_1)^{k-m-1}(\zeta-\lambda_0)^{-n-k}-\binom{m+n}{m} (\zeta-\lambda_0)^{-m-n-1}\right].
\end{eqnarray*}
On the other hand, 
$${R}_{n+1}(h^{(m)}_k, \lambda_1, \lambda_0)
=\sum_{s\in{\cal Q}_{3M_0,k} } \frac{m!}{2\pi i}\int_{\partial D_s}
h_k(\zeta) {R}_{n+1}\left(\frac{1}{(\zeta -\lambda)^{m+1}},\lambda_1, \lambda_0\right) d\zeta .
$$
Using the bound $|\zeta -\lambda_i|\geq d(\zeta, C^*)$ $(i=0,1)$ and the identity $\sum_{k=0}^m \binom{k+n-1}{k}=\binom{m+n}{m}$
we get
\begin{equation}\label{eq:rhm}
\frac{|{R}_{n+1}(h^{(m)}_k, \lambda_1, \lambda_0)|}{|\lambda_1-\lambda_0|^n}
\leq
2 H_k
\frac{m!}{2\pi }  \binom{m+n}{m} \sum_{s\in{\cal Q}_{3M_0,k} }
\int_{\partial D_s} d(\zeta, C^*)^{-m-n-1} |d\zeta|.
\end{equation}
By virtue of \eqref{eq:hm} and \eqref{eq:rhm} the proof of the theorem boils down to the following lemma
\begin{lemma}
$$
{\cal S}(k,\ell):=\sum_{s\in{\cal Q}_{3M_0,k} }
\int_{\partial D_s} d(\zeta, C^*)^{-\ell} |d\zeta| \leq C(\ell)(1+k^{2\ell}) e^k.
$$
\end{lemma}
{\bf Proof [lemma]}
First we split the sum as follows
\begin{equation}\label{split1}
{\cal S}(k,\ell)=
  \sum_{r\in{{\cal Q}^*}} \sum_{s\in V^*_r \cap{\cal Q}_{3M_0,k}}\int_{\partial D_s} d(\zeta, C^*)^{-\ell} |d\zeta|.
\end{equation}
Again
$$
d(\partial D_s,C^*) \geq \eta^{-1} d(\partial \Delta_s, K^*_d)\geq \eta_1^{-1} (1+k^{-2})^{-1/2}d( \partial V_s, K^*) \mbox{ with } \eta_1=\eta (1+k^{-2})^{1/2} .
$$
On the other hand, if $s \in V^*_r \cap {\cal Q}_{3M_0,k}$ then
$d(\partial V_s, K^*) \geq d(\partial V_{s(r,k)}, K^*)$. Since
$$\int_{\partial D_s}d(\zeta,C^*)^{-\ell} |d\zeta|\leq\eta_1^\ell d(\partial V_{s(r,k)},K^*)^{-\ell} 
\int_{\partial D_s} |d\zeta| \leq 2\eta_1^{\ell+1} d(\partial V_{s(r,k)},K^*)^{-\ell} |V_s| $$
and 
$$\sum_{s\in V^*_r\cap {\cal Q}_{3M_0,k}}|V_s|\leq |V^*_r|$$
the following estimate holds
\begin{eqnarray*}
\sum_{s\in V^*_r\cap {\cal Q}_{3M_0,k}}
\int_{\partial D_s} d(\zeta, C^*)^{-\ell}|d\zeta|
&\leq&
2 \eta_1^{\ell+1} 
d(\partial V_{s(r,k)}, K^*)^{-\ell}
\sum_{s\in V^*_r\cap {\cal Q}_{3M_0,k}} 
|V_s|\\
&\leq&
2\eta_1^{\ell+1}d(\partial V_{s(r,k)}, K^*)^{-\ell}|V^*_r| .
\end{eqnarray*}
Now we fix $n$ big enough so that it satisfies $M_{ n} \cdot \ell <1$ and $Q( n) \geq \frac{\log 8}{M_0}$.
We get the following estimates
$$
d(\partial V_{s(r,k)}, K^*) \geq
\left\{
\begin{array}{ll}
\nu_0 \frac{e^{-M_0 Q_n}}{Q^2_n} & \mbox{ if } r \in {\cal Q}^*_0(k,n),\\
\frac{e^{-M_n k}}{k^2} & \mbox{ if } r \in {\cal Q}^*_1(k,n).
\end{array}
\right.
$$
Since ${\cal Q}^*= {\cal Q}^*_0(k,n) \cup {\cal Q}^*_1(k,n)$ (see notation in the previous section)
we can split the double sum on the right hand side  of \eqref{split1} and get the estimate
\begin{eqnarray*}
{\cal S}(k,\ell) &\leq&
2 \eta_1^{\ell+1}
\left[\nu_0^{-\ell}Q^2_ne^{M_0\ell Q_n}
\left( \sum_{r \in {\cal Q}^*_0(k,n)}
|V_r|\right)
+
4 k^{2\ell}e^{M_n \ell k} \left(
\sum_{r \in {\cal Q}^*_1(k,n)}
|V_r|\right)
\right]\\
 &\leq& C_0(n,\ell)+C_1(\ell)k^{2\ell}e^k .
\end{eqnarray*}
This ends the proof of the lemma.\qed

\medskip


\section{Gevrey regularity on Diophantine points}\label{gevrey}

Let ${\cal H} $ be a Banach space and let  ${\cal G}_{\tau_0} (\lambda_0, {\cal H})$ be  the vector space of all ${\cal H}$--valued functions $h$ for which there
exist two (disjoint) open disks $\Delta^\pm$ tangent to $\Su$ at
$\lambda_0$, a formal series $\sum_{k\geq 0} c_k \Lambda^k \in
{\cal H}[[\Lambda]]$ and positive numbers $b_1$ and $b_2$ such that the
function $h$ is holomorphic in $\Delta^+\cup \Delta^-$ and it has a Gevrey-$\tau_0$ asymptotic expansion at $\lambda_0$, i.e.
\begin{equation}\label{ggrowth}
\forall N\geq 0, \ \forall \lambda \in \Delta^+\cup \Delta^-,
\|h(q)-\sum_{k=0}^{N-1}c_k (\lambda-\lambda_0)^k\|\leq b_1 b_2^N
\Gamma(1+N\tau_0)|\lambda-\lambda_0|^N,
\end{equation}
where $\Gamma$ is Euler's Gamma function.

In the following we shall slightly change our notation and set $H(\lambda,z)=H_\lambda(z)$.
It should be clear that Theorem B 
is a straightforward consequence of the following proposition:
\begin{prop}\label{pgevrey}
Let $\alpha_0 \in DC(\tau_0,\gamma_0)$ be a fixed diophantine number and let $\lambda_0=e^{2\pi i \alpha_0}$.
Let $\Delta^\pm$ be a pair of disks which
are tangent to  $\lambda_0$, $\Delta^+\subset \D$, $\Delta^-\subset \C \setminus \D$. Moreover let $\rho \in (0,e^{-1}\rho_0)$ and
${\cal H}={\cal H}^\infty(\D_\rho)$. Then
\begin{enumerate}
\item[(a)] $H\in{\cal C}^\infty_{hol}( \overline{\Delta^+\cup \Delta^-}, \  {\cal H}^\infty(\D_\rho))$  (hence $H$ is holomorphic on $\Delta^+\cup\Delta^-$).
\item[(b)] There are constants $b_1, \ b_2$ such that
\begin{equation}\label{gasintotico}
\|H - \sum_{j=0}^{N-1} \frac{1}{j!}\partial^j_\lambda H (\lambda_0,
\cdot )(\lambda-\lambda_0)^j\|_{\cal H}\leq b_1 b_2^N \Gamma(1+\tau_0 N) |\lambda-\lambda_0|^N \ \ \ \ \ \ \ \
\begin{array}{ll}\forall N \in \N, \\
 \forall \lambda \in \Delta^\pm.
\end{array}
\end{equation}
\end{enumerate}
\end{prop}

\subsection{More Cantor sets related with the diophantine condition}
In order to prove Proposition \ref{pgevrey} we shall use just a few
definitions and results from the paper \cite{MS1}\footnote{Let us point out that we do not stick to the
notation used in \cite{MS1}: in particular we shall call
``$K_{\gamma,\tau}$'' the set that in \cite{MS1} is called
$C_{\psi_{\gamma, \tau}}$.}. For $\tau\geq 2$, $\gamma \in (0,1)$ let
$$K_{\gamma,\tau}:= \{ x \in (\R \setminus \Q)/\Z \ : \ \forall k\geq 0 \ q_{k+1}\leq \gamma^{-1}q_k^{\tau-1} \};$$
It is readily seen that $K_{\gamma,\tau}$ is a compact subset of $\R/\Z$ and
\begin{equation}\label{diofanto}
D(\gamma, \tau)\subset K_{\gamma,\tau} \subset D(\frac{\gamma}{2} , \tau) \mbox{ }
\end{equation}
(see also A3.2 in \cite{MS1}).
The proof of the following proposition can be found in  \cite{MS1} (Proposition 2.2).
\begin{prop}\label{DV}.\\
\begin{enumerate}
\item Each connected component $]\xi^-, \xi^+[$ of $(\R/\Z)\setminus
 K_{\gamma,\tau}$ contains a unique rational number $p/q$
 which is a convergent of both endpoints $\xi^\pm$. We shall call such
 convergent $p/q$ the {\em pseudocenter} of the component $]\xi^-,
 \xi^+[$. We shall denote ${\cal Q}_{\gamma, \tau}$ the set of
 pseudocenters of all connected components of $(\R/\Z)\setminus
 K_{\gamma,\tau}$.
\item
$\frac{\gamma}{2}q^{-\tau} \leq |\xi^\pm -\frac{p}{q}| \leq 2\gamma q^{-\tau}$.
\end{enumerate}
\end{prop}

\medskip
Let us remark that, if $\alpha \in K_{\gamma,\tau}$ then
$q_{k+1}=a_{k+1}q_k+q_{k-1}\leq \gamma^{-1} q^{\tau -1}$ and hence
$a_{k+1}\leq \gamma^{-1} q^{\tau -2}$. From this information we get
not only an a priori estimates for ${\cal B}(\alpha)$ but also
$$\sum_{k=n}^{+\infty}\frac{\log a_{k+1}}{q_k}\leq (\log \gamma^{-1})\sum_{k=n}^{+\infty}\frac{1}{q_k}+
(\tau -2)
\sum_{k=n}^{+\infty}\frac{\log q_{k}}{q_k}.$$
This is interesting because, since $q_k \geq F_k$, we get that
$$\sum_{k=n}^{+\infty}\frac{1}{q_k}\leq \sum_{k=n}^{+\infty}\frac{1}{F_k}
,  \ \ \ \ \ \ \ \sum_{k=n}^{+\infty}\frac{\log q_{k}}{q_k} \leq \sum_{k=n}^{+\infty}\frac{\log F_{k}}{F_k};$$
therefore setting
\begin{equation}\label{emen}
M_n:=(\log \gamma^{-1})\sum_{k=n}^{+\infty}\frac{1}{F_k}+
(\tau -2)
\sum_{k=n}^{+\infty}\frac{\log F_{k}}{F_k}
\end{equation}
 we see that $K_{\gamma, \tau} \subset K^*_{(M_n)}:= \{ x\in \R/\Z \ :\ \sum_{k=n}^{+\infty} \frac{\log a_{k+1}}{q_k}\leq M_n\}$.
  For $s \in {\cal
Q}_{\gamma, \tau}$ let $V_s$ be the corresponding connected component
of $(\R/\Z)\setminus K_{\gamma, \tau}$, $\Delta_s$ be the diamond over
$V_s$, $D_s:=\exp^\#(\Delta_s)$ and $C_{\gamma, \tau} = \C \setminus
\cup_{s\in {\cal Q}_{\gamma,\tau}} D_s$. We then have that $C_{\gamma,
\tau }$ is contained in the set of ${\cal C}^\infty$-regularity
$C^*_{(M_n)}$.

The following geometrical lemma is useful to settle both {\em (a)} and {\em (b)} of Proposition \ref{pgevrey}.

\begin{lemma}\label{gl}
Let $\tau_0 \geq 2$, $\alpha_0 \in DC(\tau_0,\gamma_0)$,   $\lambda_0=e^{2\pi i \alpha_0}$. Let $Q$ be a closed set satisfying
\begin{enumerate}
\item[(i)] $Q\cap \Su = \{\lambda_0\}$;
\item[(ii)]
 there exist $\mu_0>0 $ such that ${\rm dist}(\lambda, \Su)\geq \mu_0|\lambda -\lambda_0|^2 \ \ \forall \lambda \in Q$.
\end{enumerate}
Then for any fixed $\tau>2\tau_0$ there exists $\gamma \in (0,\gamma_0)$ such that,
defining $(M_n)$ as in \eqref{emen} so that $K_{\gamma, \tau}\subset
C^*_{(M_n)}$,
\begin{enumerate}
\item[(a')] $Q\subset K_{\gamma, \tau}\subset C^*_{(M_n)}$;
\item[(b')]
 there exists  $\mu>0 $ such that ${\rm dist}(\zeta, Q)\geq \mu \gamma_0^{2} q^{-2\tau_0} \ \ \forall  \zeta \in \partial D_{p/q}, (p/q \in {\cal Q}_{\gamma, \tau})$.
\end{enumerate}
\end{lemma}
{\bf Proof [lemma]} Let us consider $Q_0:=Q \cap \{z:
|z-\lambda_0|\leq 1/2 \}$ and $Q_1:=Q \cap \{z: |z-\lambda_0|\geq 1/2 \}$,
we shall prove that (a) and (b) hold on both the closed sets $Q_0$ and
$Q_1$ and hence hold on $Q$ as well.  Of course, in the case of $Q_1$
there is no problem: since $Q_1\cap \Su = \emptyset$ the points of
$Q_1$ are bounded away from $\Su$ and both (a) and (b) hold provided that
 $\gamma$ and $\mu$ are  small enough.
As far as $Q_0$ is concerned, we observe that the logarithm is well defined on $\{z:
|z-\lambda_0|\leq 1/2 \}$ and it has a bounded distorsion property, therefore wecan check both (a) and (b) for $Q_0$ just proving
  the following statement:
\begin{quote}
If $\mu_0'> 0$ and $Q':=\{ \alpha : \Im \alpha \geq\mu_0'
  |z-\alpha_0|^2, \ |\Re \alpha |\leq \pi/6\}$ there exist $\gamma,
  \mu' >0$ such that ${\rm dist}(\xi, Q')\geq \mu' \gamma_0^2
  q^{-2\tau_0} \ \ \forall \xi \in \partial \Delta_{p/q}, p/q \in
  {\cal Q_{\gamma, \tau}}$.
\end{quote}
This is readily checked because, if $\xi \in \partial \Delta_{p/q}$ is
fixed and $\alpha(\xi)$ is the nearest point in $Q'$ we have that
$|\alpha(\xi)-\xi|\geq |\alpha(\xi)-p/q|-|\xi - p/q|$; on the other
hand we have that
$$|\alpha(\xi)-p/q| \geq C|\Im \alpha|\geq
C\mu_0'|\alpha -p/q|^2\geq (C\mu_0' \gamma_0^2)q^{-2\tau_0}$$ while
$$
|\xi -p/q|\leq 2\frac{\kappa}{\sqrt{1+\kappa^2}}\gamma q^{-\tau} = O(q^{-\tau})
$$
\qed

\medskip

Let us point out that by means of the same argument used in the proof of the last lemma we also get
$$|\lambda_0 - \zeta| \geq \mu \gamma_0 q^{-\tau_0}\ \ \forall \zeta \in
\partial D_{p/q}, (p/q \in {\cal Q_{\gamma, \tau}}).$$

\subsection{Proof of Gevrey regularity}
By (a') of Lemma \ref{gl} we get that, if $\tau >2\tau_0$, it is
possible to choose $\gamma$ such that $\overline{\Delta^+\cup
\Delta^-}\subset K_{\gamma, \tau}\subset C^*_{(M_n)}$ hence $H$ is
${\cal C}^\infty_{hol}(\overline{\Delta^+\cup \Delta^-}, {\cal H}^\infty(\D_\rho))$ in
particular $H$ admits a Taylor expansion at $\lambda_0$.
We shall set 
$${\bf R}_N(H,\lambda, \lambda_0):= H(\lambda, \cdot) - \sum_{j=0}^{N-1}
\frac{1}{j!}\partial^j_\lambda H (\lambda_0, \cdot )(\lambda-\lambda_0)^j$$ so
that ${\bf R}_N(H, \cdot, \lambda_0):\Delta^+\cup \Delta^- \to
{\cal H}^\infty(\D_\rho)$.  We point out that, setting
$${R}_{N}(h_k, \lambda, \lambda_0):= h_k(\lambda)- \sum_{j=0}^{N-1}
 \frac{1}{j!}h_k^{(j)} (\lambda_0)(\lambda-\lambda_0)^j,$$
we can write
\begin{equation}\label{modes}
{\bf R}_N(H,\lambda, \lambda_0)(z)=\sum_{k=1}^{+\infty}{R}_N(h_k, \lambda, \lambda_0) z^k
\end{equation}

\bigskip

Let us recall that, by Cauchy formula,
$${R}_N(h_k, \lambda, \lambda_0)=\sum_{s\in{\cal Q}_{\gamma, \tau, k} } \frac{1}{2\pi i}\int_{\partial D_s}
 \frac{h_k(\zeta)}{(\zeta -\lambda_0)^{N}(\zeta -\lambda)}(\lambda - \lambda_0)^N d\zeta,
$$ where ${\cal Q}_{\gamma, \tau, k}={\cal Q}_{\gamma, \tau}\cap
\Q_k$.  Moreover, if $\lambda \in \Delta^+\cup \Delta^-$ and $\zeta
\in \partial D_s$ with $s=p/q \in{\cal Q}_{\gamma, \tau,k}$ then 
$$|h_k(\zeta)|\leq
H_k, \ \ \ \ \ \ \ |\zeta -\lambda_0|\geq \mu\gamma_0|q|^{-\tau_0} \ \ \ \ {\rm and } \ \ \ \ \  |\zeta -
\lambda|\geq \mu\gamma_0^2 |q|^{-2\tau_0}$$
 so we get the bound
$$|{R}_N(h_k,\lambda, \lambda_0)|\leq
(\mu \gamma_0)^{-N-2}
k^{(N+2)\tau_0}
H_k
\frac{1}{2\pi} \sum_{s\in{\cal Q}_{3M_0,k} }
\int_{\partial D_s} |d\zeta| |\lambda -\lambda_0|^N$$
The sum on the right hand side of the last formula is bounded by a constant (independent from $k$ and $N$)
moreover, since we have chosen $\rho<\rho_0$, we can fix $C_1$ such that
 $k^{2\tau_0} H_k \rho^k \leq C_1 e^{-k}$.
Thus, summing up on $k$, we get
$$\|{\bf R}_N(H,\lambda, \lambda_0)\|\leq \sum_{k=1}^{+\infty}|{R}_n(h_k, \lambda, \lambda_0)| \rho^k \leq C_2
(\mu\gamma_0)^{-N}\left(\sum_{k=0}^{+\infty} k^{N\tau_0} e^{-k}\right) |\lambda -\lambda_0|^N$$
where $C_2$ is a suitable constant. Since
$$
\sum_{k=0}^{+\infty} k^{\beta}e^{-k}\leq \int_0^{+\infty}t^\beta e^{-t}dt +
(\beta/e)^\beta \leq 2\Gamma(1+\beta)
$$
the thesis follows. \qed

\section{ Appendix A. Arithmetical tools and Brjuno series}

\vskip .5 truecm
Let us recall some notation and elementary facts about classical continued fractions. We refer the reader to \cite{HW} and \cite{Khi} for more details.
Continued fractions  are obtained by coding the orbits of real numbers under the
iteration of the Gauss map
${\cal G} : (0,1) \mapsto [0,1]$
defined by
${\cal G} (x) = \{x^{-1}\} = x^{-1} -
[x^{-1}]$ where $[x] $ and $\{ x\}$ respectively denote the integer and the fractional part of $x$.
This map is piecewise analytic with inverse branches $T_n(x)=\frac{1}{n+x}$,
$T_n={\cal G}^{-1}$ on the interval $\left(\frac{1}{n+1},\frac{1}{n}\right)$.
Given $x \in \mathbb R \setminus \mathbb Q$ we set
$x_0 = x - [x] \; ,            a_0  = [x] \; ,$
then one obviously has
$x = a_0 +  x_0$.
We now define inductively for all $n \ge 0$
$
x_{n+1} = {\cal G}(x_n)$,
$a_{n+1} = [ x_n^{-1}]  \ge 1$,
thus
$x_{n} = T_{a_{n+1}}(x_{n+1})$.
Therefore we have
$$
x=a_0+T_{a_1}(x_1)=\ldots =a_0+T_{a_1}\circ ... \circ T_{a_n}(x_n)
     =a_0 + \displaystyle\frac{1}{ a_1
     + \displaystyle\frac{1  }{ a_2 + \ddots +
     \displaystyle\frac{1 }{ a_n + x_n}}}\; .
$$

We will use the short notation
$
x=[a_0,a_1,\ldots ,a_n,
     \ldots]
$ for the infinite fraction.
The nth-convergent is then the rational number corresponding to the
finite fraction
$\frac{p_n }{ q_n} = [a_0,a_1,\ldots ,
                     a_n] $.

The numerators $p_n$ and denominators
$q_n$ are recursively determined
for all $n \ge 0$ by
$$
p_n = a_n p_{n-1} +  p_{n-2} \; ,
            q_n = a_n q_{n-1} + q_{n-2} \; ,
$$
with the initial conditions
$
p_{-1}=q_{-2}=1 \;\;,\;\;\;p_{-2}=q_{-1}=0
$. Note that $ q_n p_{n-1} - p_n q_{n-1} = (-1)^n $.

For all $n\ge 0$ one also has
$$
x =\frac {p_n + p_{n-1}  x_n }{q_n + q_{n-1}
       x_n } \; ,  \;\;\;\;
x_n = - \frac{q_n x -p_n}{ q_{n-1} x - p_{n-1}}\; ,
$$
thus
for all $k\ge 0$ and for all $x\in\R\setminus\Q$ one has
$\frac{p_{2k}}{ q_{2k}}<x<\frac{p_{2k+1}}{q_{2k+1}}$.

It is not difficult to show that for all $x \in \mathbb R \setminus \mathbb Q$ and for all $n \ge 1$ one has
$q_n\ge\frac{1}{2}\phi_0^{1-n}$, with
$
\phi_0= \frac{\sqrt{5}-1}{ 2}
$.
This implies that the series
$\sum_{k=0}^\infty \frac{\log q_{k}}{ q_{k}}$ and
$\sum_{k=0}^\infty \frac{1}{ q_{k}}$
are always convergent and that their sum is uniformly bounded.

\vskip .3 truecm
\noindent
For all integers  $k\ge 1$, the iteration of the Gauss map $k$ times leads
to the following partition of
$(0,1)$; $\sqcup_{a_{1},\ldots ,a_{k}}
I(0,a_{1},\ldots ,a_{k})$, where $a_{i}\in {\mathbb N}$,
$i=1,\ldots ,k$, and
\begin{equation}\label{cylinders}
I(0,a_{1},\ldots ,a_{2k}) =
\left( \frac{p_{2k}}{ q_{2k}},\frac{p_{2k}+p_{2k-1}}{ q_{2k}+q_{2k-1}}\right)
\; , \;\;\;
I(0,a_{1},\ldots ,a_{2k+1}) =\left( \frac{p_{2k+1}+p_{2k}}{ q_{2k+1}+q_{2k}}, \frac{p_{2k+1}}{ q_{2k+1}}\right).
\end{equation}
These are the intervals corresponding to the branches of ${\cal G}^k$: they are determined by the fact that all points
$x\in I(0,a_{1},\ldots ,a_{k}) $ have the
first $k+1$ partial quotients exactly equal to $\{0,a_{1},\ldots ,a_{k}\}$.
Thus
$$
I(0,a_{1},\ldots ,a_{k}) = \left\{ x\in (0,1)\,\mid\; x=
\frac{p_{k}+p_{k-1}y}{ q_{k}+q_{k-1}y}\; , \; y\in (0,1)\right\}\; .
$$
Note that $\frac{dx}{ dy}= \frac{(-1)^{k}}{ (q_{k}+q_{k-1}y)^{2}}$ is
positive (negative) if $k$ is even (odd). It is immediate to check
that any rational number $p/q\in (0,1)$, $(p,q)=1$, is the endpoint
of exactly two branches of the iterated Gauss map. Indeed $p/q$ can
be written as $p/q=[0,\bar{a}_{1},\ldots ,\bar{a}_{k}]$ with $k\ge 1$
and $\bar{a}_{k}\ge 2$ in a unique way and it is the left (right)
endpoint of $I(0,\bar{a}_{1},\ldots ,\bar{a}_{k})$ and
the right (left) endpoint of $I(0,\bar{a}_{1},\ldots ,\bar{a}_{k}-1,1)$
if $k$ is even (odd).

The cylinders $I(a_0,...,a_N)=a_0+I(0,a_1,...,a_N)$ form a partition
of the whole real line as $a_0$ varies in $\mathbb Z$ and $(a_1,\ldots ,a_n)\in
\mathbb N^n$.

\section{ Appendix B. ${\cal C}^{1}$-holomorphic and ${\cal C}^{\infty}$-holomorphic functions}

\vskip .3 truecm
Let $(B,\Vert\;\Vert )$ be a complex Banach space. In this appendix
we recall the definition of ${\cal C}^{1}$-holomorphic and ${\cal C}^{\infty}$-holomorphic functions
as they are given respectively in \cite{He} and \cite{Ris}. We follow quite closely Section 2 of \cite{MS1}
to which we refer for a more detailed discussion.

\vskip .3 truecm
Let $C$ be a  compact subset of $\C$ or of~${\mathbb P}^1\C$.
If $C\subset \C$, a continuous function $f:\,C\rightarrow B$ is said to be ${\cal C}^{1}$-holomorphic
if there exists a continuous map $f^{(1)} :\, C\rightarrow B$ such that
\begin{eqnarray*}
\forall \lambda\in C,\;\forall\varepsilon >0, \;
 \exists\delta >0 \;/\quad
\forall \lambda_{1},\lambda_{2}\in C,\;
&|\lambda_1-\lambda|<\delta,\; |\lambda_2-\lambda|<\delta
\\
&\qquad\Rightarrow
\Vert f(\lambda_{2})-f(\lambda_{1})-f^{(1)}(\lambda_{1})(\lambda_{2}-\lambda_{1})
\Vert \le \varepsilon |\lambda_{1}-\lambda_{2}|.
\end{eqnarray*}
This definition extends in an obvious way to the case $C\subset {\mathbb P}^1\C$ by means of
the standard complex coordinates charts.

The above definition makes use of
the generalization of the notion of smoothness of a function to
a closed set due to Whitney \cite{Wh}. Notice however
 that $f^{(1)}$  is a complex derivative:
$\bar{\partial}f =0$,
$\partial f = f^{(1)}$
and $f$ is holomorphic in the interior of~$C$.

The space ${\cal C}^{1}_{hol}(C,B)$
becomes a  Banach space by taking as norm
$$
||| f |||= \max\Bigl(\sup_{\lambda\in C}\Vert f(\lambda)\Vert\, , \,
\sup_{\lambda\in C}\Vert f^{(1)}(\lambda)\Vert\, , \,
\sup_{\lambda_{1},\lambda_{2}\in C, \, \lambda_{1}\not=\lambda_{2}}\frac{\|
f(\lambda_{2})-f(\lambda_{1})-f^{(1)}(\lambda_{1})(\lambda_{2}-\lambda_{1})\|}{|\lambda_{1}-\lambda_{2}|}
\Bigr)
$$

Let ${\cal R}(C,B)$
denote the uniform algebra of continuous functions
from $C$ to $B$ which are uniformly approximated
by rational functions with all the poles outside $C$. A very important
property of the space ${\cal C}^{1}_{hol}(C,B)$ is that it is a linear subspace of
${\cal R}(C,B)$. This fact allows to prove that functions in
${\cal C}^{1}_{hol}(C,B)$
share some of the properties of holomorphic functions.
If $(U_{\ell})_{\ell\ge 1}$ denote the connected components of
${{\mathbb P}^1\C}\setminus C $, assuming that each $\partial U_{\ell}$
is a piecewise smooth Jordan curve and $\sum_{\ell\ge 1}{\rm length}
(\partial U_\ell)<+\infty$, Cauchy's theorem holds:
$$
\sum_{\ell=1}^\infty \int_{\partial U_\ell} f(\lambda)\,d\lambda=0.
$$
This is very easy to see: since $f\in {\cal R}(C,B)$,
one can approximate $f$ by a sequence $(r_{k})_{k\in N}$
of $B$-valued
rational functions with poles off $C$. Cauchy's theorem applies to
these rational functions and one can pass to the limit because the
convergence is uniform.
Moreover, at all points $\lambda\in C$ such that
$$
\sum_{\ell=1}^\infty \int_{\partial U_\ell} \frac{|d\zeta|}{ |\zeta-\lambda|}<+\infty,
$$
Cauchy's formula also holds:
$$
f(\lambda) = \frac{1}{ 2\pi i} \sum_{\ell=1}^\infty \int_{\partial U_\ell}
\frac{f(\zeta)}{\zeta-\lambda}\,d\zeta.
$$

One can also
define higher order derivatives by means of Cauchy's formula, but in order to do so
one needs further assumptions on~$\lambda$
(namely $\sum_{\ell=1}^\infty \int_{\partial U_\ell} \frac{|d\zeta|}{
|\zeta-\lambda|^{n+1}}<+\infty $ to obtain a derivative of order $n$).

\vskip .3 truecm

A function $f:\,C\rightarrow B$ is said to be
${\cal C}^{\infty}$-holomorphic
if there exist an infinite  sequence of continuous functions
$(f^{(n)})_{n\in {\N}}:\,C\rightarrow B$ (the ``$n$-th complex
derivative of $f$'')
such that $f^{(0)}= f$
and, for all $n,m\ge0$, the function~$R^{(n,m)}$ defined by
$$
 R^{(n,m)}(\lambda_{1},\lambda_{2}) =
\sum_{h=0}^{m}
\frac{f^{(n+h)}(\lambda_1)}{h!}(\lambda_{2}-\lambda_{1})^h -
f^{(n)}(\lambda_{2}), \qquad \lambda_{1},\lambda_{2}\in C, $$
has the following property:
$$
\forall \lambda\in C,\;\forall\varepsilon >0,\; \exists\delta >0 \;/\;
\forall \lambda_{1},\lambda_{2}\in C,\;
|\lambda_1-\lambda|<\delta,\; |\lambda_2-\lambda|<\delta
\;\Rightarrow\;
\Vert R^{(n,m)}(\lambda_1,\lambda_2)\Vert\le \varepsilon |\lambda_1-\lambda_2|^{m}.
$$

The space of
${\cal C}^{\infty}$-holomorphic $B$-valued functions
on a compact set is a Fr\'echet space.
We stress once more that the derivatives are taken in a complex sense,
thus $\bar{\partial} f^{(n)}=0$ for all $n\in {\N}$.
The functions~$f^{(n)}$ are some generalized ``weak complex derivatives for
$f$''; clearly
$f$ must be analytic in the interior of $C$ and
$$
\forall n,m\in {\N},\quad
\forall \lambda\in {\rm int}(C),
\quad
f^{(n+m)}(\lambda) = \partial^m f^{(n)}(\lambda).
$$

\vskip .5 truecm

Let $(C_{j})_{j\in {\N}}$ be a monotonic non-decreasing sequence of compact
subsets of~${{\mathbb P}^1\C}$. The associated space of {\sl $B$-valued monogenic functions}
is defined to be the projective limit
$$
{\cal M}((C_{j}),B) = \varprojlim {\cal C}^{1}_{hol}(C_{j},B).
$$

The restrictions ${\cal C}^{1}_{hol}(C_{j+1},B)\rightarrow
{\cal C}^1_{hol}(C_{j},B)$ are continuous linear operators between Banach spaces,
thus ${\cal M}((C_{j}),B)$ is a Fr\'echet space with seminorms
$\Vert \,.\, \Vert_{{\cal C}^{1}_{hol}(C_{j},B)}$.

The above definition is inspired by the work of Borel \cite{Bo}
(see also \cite{He}, p.\  81). Borel considered the case $B={\C}$ and wanted to extend
the notions of holomorphic function and analytic continuation.
Borel's idea was to allow monogenic continuation through natural boundaries of
analyticity
by selecting points at which
the function is ${\cal C}_{hol}^1$-holomorphic. If the function is moreover
${\cal C}_{hol}^\infty$-holomorphic at such a point, the question of quasianalyticity may
be raised: Is the function determined by its Taylor series?
Such a uniqueness property could depend on the choice of the sequence~$(C_j)$ which defines the
monogenic class (and not only on the union of the~$C_j$'s),
and the Cauchy formula could help to establish it.

Unfortunately this strong form of quasianalyticity is not true in
general spaces of monogenic functions unless some rather restrictive
assumptions are made (see, e.g.\ , \cite{Wk}). However it is proved in
\cite{MS2} that it is quite a general property that the spaces of
${\cal C}^1$--holomorphic functions (and also of monogenic functions)
which appear in linearization problems have a weaker quasianalyticity
property, namely their functions cannot vanish on a set of positive
$1$--dimensional Hausdorff measure without being identically equal to
zero. This provides an example of generalized analytic continuation
(for a comprehensive discussion of generalized analytic continuations
{\it other} than Borel's theory see \cite{RS}).

The notion of $B$-valued monogenic function is well adapted to
linear small denominator problems, as the cohomological equation
considered in \cite{MS1} but it is useless in nonlinear problems
since one cannot fix a target Banach space if one wants the
increasing sequence of compact sets to include all points on the
unit circle verifying the Brjuno condition. Indeed the radius of
convergence of the linearization is, in general, also bounded above in terms of
the exponential of minus the Brjuno function as first proved by
Yoccoz \cite{Yo1}. For this reason we introduce the notion of {\it
monogenic function with values in} $\C\{z\}$: suppose
$(C_j)_{j\in\N}$ is a monotonic non-decreasing sequence of compact
subsets of~${\mathbb P}^1\C$ and consider the increasing sequence
of Banach spaces ${\cal H}^\infty (\D_{r_j})$ associated to a
monotonic non-increasing sequence of radii $r_j \to 0$. Consider
the Banach space $B_\ell=\cap_{j=0}^\ell {\cal C}^1_{hol} (C_j, {\cal
H}^\infty (\D_{r_j}))$ with the norm $\Vert f\Vert_{B_\ell} =
\max_{0\le j\le \ell} \Vert f\Vert_{{\cal C}^1_{hol} (C_j, {\cal
H}^\infty (\D_{r_j}))}$. Clearly the injections $i_\ell\, :\,
B_\ell\hookrightarrow B_{\ell-1}$ are bounded linear operators between
Banach spaces with norms $\Vert i_\ell\Vert\le 1$. The projective
limit of the system of Banach spaces
$$ {\cal M}( (C_j),
\C\{z\}) = \varprojlim B_j
$$
is the space of monogenic functions with values in the holomorphic
germs $\C\{z\}$. It is a Fr\'echet space with the seminorms
$\Vert \cdot\Vert_{B_\ell}$.
Thus, as a set, ${\cal M}( (C_j),
\C\{z\})$ consists of all
the functions which are defined in
$C = \bigcup_{j\in\N} C_j$
and such that, for every $j\in\N$, the restriction $f_{|C_j}$
belongs to ${\cal C}^1_{hol} (C_j, {\cal H}^\infty (\D_{r_j})).$
This space, being the projective limit of the Banach spaces $B_\ell$, may depend on the increasing sequence~$(C_j)$ and on the decreasing sequnce $r_j$ 
(rather than
on the set~$C$ only).

\end{document}